\documentclass[12pt]{amsart}
\usepackage[dvips]{graphics}
\usepackage{amsfonts}
\usepackage{amssymb}
\usepackage{a4}

\DeclareMathAlphabet{\eusm}{U}{}{}{}  
\SetMathAlphabet\eusm{normal}{U}{eus}{m}{n}
\SetMathAlphabet\eusm{bold}{U}{eus}{b}{n}

\DeclareMathAlphabet{\eufrak}{U}{}{}{}  
\SetMathAlphabet\eufrak{normal}{U}{euf}{m}{n}
\SetMathAlphabet\eufrak{bold}{U}{euf}{b}{n}

\newtheorem{theorem}{Theorem}[section]
\newtheorem{proposition}[theorem]{Proposition}
\newtheorem{lemma}[theorem]{Lemma}
\newtheorem{corollary}[theorem]{Corollary}

\theoremstyle{definition}
\newtheorem{definition}[theorem]{Definition}

\theoremstyle{remark}
\newtheorem{remark}[theorem]{Remark}
\numberwithin{equation}{section}

\newcommand{\bmultg}{\!\begin{array}{c} {\scriptstyle\times} \\[-12pt]\cup\end{array}\!}
\newcommand{\bmultsqrt}{\!\begin{array}{c} {\scriptstyle\times} \\[-12pt]\tilde{\cup}\end{array}\!}
\newcommand{\bmultk}{\!\!{\scriptstyle\begin{array}{c} {\scriptscriptstyle\times} \\[-12pt]\cup\end{array}}\!\!}

\newcommand{\mmultalt}{\kern0.3em\rule{0.04em}{0.52em}\kern-.15em\tilde{\gtrdot}\kern0.2em}

\newcommand{\mmultg}{\kern0.3em\rule{0.04em}{0.52em}\kern-.35em\gtrdot}
\def\mmultk{\kern0.17em\lower0.1ex\hbox{\rule{0.025em}{0.43em}}\kern-.105em\gtrdot}

\begin{document}

\title[Monotone and Boolean Convolutions]{Monotone and Boolean Convolutions for Non-compactly Supported Probability Measures}
\author{Uwe Franz}
\address{D\'epartement de math\'ematiques de Besan\c{c}on,
Universit\'e de Franche-Comt\'e 16, route de Gray, 25 030
Besan\c{c}on cedex, France}
\curraddr{Graduate School of Information
Sciences, Tohoku University, Sendai 980-8579, Japan}
\urladdr{http://www-math.univ-fcomte.fr/pp\underline{ }Annu/UFRANZ/}
\thanks{This work was supported by a Marie Curie Outgoing International
Fellowship of the EU (Contract Q-MALL MOIF-CT-2006-022137), an ANR Project (Number ANR-06-BLAN-0015), and a
Polonium cooperation}

\begin{abstract}
The equivalence of the characteristic function approach and the probabilistic approach to monotone and boolean convolutions is proven for non-compactly supported probability measures. A probabilistically motivated definition of the multiplicative boolean convolution of probability measures on the positive half-line is proposed. Unlike Bercovici's multiplicative boolean convolution it is always defined, but it turns out to be neither commutative nor associative. Finally some relations between free, monotone, and boolean convolutions are discussed.
\end{abstract}

\maketitle

\section{Introduction}
There are at least three equivalent approaches to (classical) convolutions of probability measures on groups and semi-groups, which we shall call here the harmonic analysis, the probabilistic, and the characteristic functions approach.  The main goal of the present paper is to prove the equivalence of the probabilistic and the characteristic functions approach for the convolutions based on monotone and boolean independence introduced recently.

In harmonic analysis the convolution of two probability measure $\mu$ and $\nu$ on a group or semi-group $G$ is defined as the image measure $T^{-1}\mu\otimes\nu$ of their product $\mu\otimes\nu$ under the map $T:G\times G\to G$ defined by the (semi-) group operation, $T(x,y)=xy$. The probabilistic definition would be to take the law of the product $XY$ of two $G$-valued random variables $X$ and $Y$ as the convolution of $\mu$ and $\nu$, if $X$ and $Y$ are independent and distributed according to $\mu$ and $\nu$, respectively. The equivalence of the two approaches follows from the fact that two random variables are independent if and only if their joint law is the product of their marginals.

In the characteristic function approach one defines a bijection between probability measures and an appropriate class of functions, e.g.\ in the case of probability measures on the real line one usually chooses the Fourier transform, $\hat{\mu}(t)=\int_\mathbb{R} e^{itx}{\rm d}\mu(x)$. By Bochner's theorem, the Fourier transform defines a bijection between the probability measures on the real line and the class of uniformly continuous, positive definite functions one the real line, whose value at the origin is equal to one. Checking that this class is closed under the pointwise products, one could define the (additive) convolution of the probability measures $\mu$ and $\nu$ on $\mathbb{R}$ as the unique probability measure $\lambda$ such that $\hat{\lambda}(t) = \hat{\mu}(t)\hat{\nu}(t)$ for all $t\in\mathbb{R}$. Usually this approach is not used to define the convolution but rather to compute it or to study its properties, prove limit theorems, etc.

These three approaches can also be used to define and study the free, monotone, and boolean convolutions, which are based on the corresponding notions of independence in quantum probability. For this purpose groups have to be replaced by dual groups \cite{voiculescu87,voiculescu90} and random variables have to be replaced by quantum random variables, i.e.\ appropriate classes of operators on Hilbert spaces. Furthermore, it is necessary to find adequate characteristic functions.

In the free case the equivalence of the probabilistic and the characteristic function approach was first proven in \cite{voiculescu86} for the additive free convolution on the real line and for probability measures with compact support. Later it was extended to multiplicative free convolutions \cite{bercovici+voiculescu92} and probability measures with non-compact support in \cite{maassen92,bercovici+voiculescu93}.

For the convolutions based on the monotone and boolean independence, this equivalence has been proven so far only for probability measures with compact support, cf.\ \cite{muraki00,bercovici04,bercovici05,franz05a} and \cite{speicher+woroudi93,franz04,bercovici04b}. For non-compactly supported probability measures these convolutions are currently defined only via their characteristic functions (which are certain functions of their Cauchy transforms), and not via the more natural harmonic analysis or probabilistic approach. In this paper, we will extend the equivalence between the probabilistic approach and the characteristic function approach to non-compactly supported probability measures.

This extension is necessary, e.g., for studying monotone and boolean counterparts of stable laws by probabilistic means. We also need it to extend the results of \cite{franz+muraki05} beyond the bounded case, and to construct and study unbounded quantum stochastic processes with monotonically or boolean independent increments. In the monotone case, these processes turn out to have an interesting relation to the theory of Loewner chains, cf.\ \cite{franz07a}.

The equivalence between the harmonic analysis approach and the probabilistic approach for the free convolutions was established in \cite{voiculescu87,voiculescu90} where dual groups where introduced for this purpose. Due to the universal properties of the products underlying the free, monotone, and boolean independence, it follows in the same way for the monotone and boolean convolutions, cf. \cite{schuermann95b,benghorbal+schuermann02,franz03c,muraki03}.

To apply the probabilistic approach to non-compactly supported probability measures, it is necessary to extend the definition of monotone or boolean independence to unbounded operators. This is done in Definitions \ref{def-mon-indep} and \ref{def-bool-indep}. Two not necessarily bounded operators that admit a functional calculus for continuous functions are called monotonically or boolean independent, if the algebras consisting of bounded functions vanishing at the origin of these operators are monotonically or boolean independent. The restriction to functions vanishing at the origin is necessary, because for monotone and boolean independence one has to allow non-unital algebras to get non-trivial examples, see also \cite{benghorbal+schuermann02,muraki03}.

As a next step we show that the general case of two monotonically or boolean independent normal operators can be reduced to a universal model which allows explicit calculations, see Theorems \ref{thm-mon-model} and \ref{thm-bool-model}.

In Theorems \ref{thm-mon-add-conv} and \ref{thm-bool-add-conv} we then show that the sum of two monotonically or boolean independent self-adjoint operators is essentially self-adjoint, if the state vector is cyclic for the algebra generated by these two operators. This condition is not very restrictive. If the state vector is not cyclic, then one can always restrict to the subspace that is generated from it. Furthermore, we show that the distribution of the sum is equal to the additive monotone or boolean convolution of the distributions of the two operators. This justifies the name {\em additive} monotone or boolean convolution.

Next we treat the multiplicative convolutions of probability measures on the positive half-line. This case is more subtle, because there are many possibilities for constructing ``multiplicatively'' a positive operator out of two given positive operators. Here we consider only the two cases $(X,Y)\mapsto \sqrt{X}Y\sqrt{X}$ and $(X,Y)\mapsto \sqrt{Y}X\sqrt{Y}$. These operations  are neither commutative nor associative. Nonetheless, in the free case, both lead to the same multiplicative free convolution of probability measures on the positive half-line, which is associative and commutative. This follows from the fact that the free product is trace-preserving.

But the monotone and the boolean product are not trace-preserving, and so it is not surprising that the situation becomes more complicated here.

In the monotone case, we show that the operation $(X,Y)\mapsto \sqrt{X}Y\sqrt{X}$ leads to a ``nice'' convolution product that agrees with the one defined in \cite{bercovici04}, see also \cite{franz05a} and Theorem \ref{thm-mon-mult}.

We also study the case $(X,Y)\mapsto \sqrt{Y}X\sqrt{Y}$, but it turns out that this leads to a convolution which is not associative and which does not seem to have a ``nice'' definition in terms of characteristic functions, see Subsection \ref{subsec-mon-mult-other}.

Then we consider these operations for boolean independent positive operators and give an alternative definition of the multiplicative boolean convolution for probability measures on the positive half-line. The definition proposed by Bercovici \cite{bercovici04b} has the disadvantage that it is not always defined. In Remark \ref{rem-new-bool-mult-conv}, we give a new probabilistically motivated definition of the multiplicative convolution of two probability measures on the positive half-line. But it turns out that this convolution is neither commutative nor associative.

Probability measures on the unit circle are of course always compactly supported and the associated operators are unitary and therefore bounded. In this case the equivalence between the probabilistic and the characteristic function approach was already established in \cite{bercovici04,franz05a,franz04}. For completeness we recall these results in Subsections \ref{subsec-mult-mon-T} and \ref{subsec-mult-bool-T}.

Finally, in Section \ref{sec-relations} we discuss some relations between free, monotone, and boolean convolutions.

\section{Preliminaries}

\subsection{Quantum probability}

In quantum probability the commutative algebra of random variables is replaced by a possibly non-commutative algebra $\mathcal{A}$ of operators on a Hilbert space $H$ and the role of the probability measure is taken by a state $\Phi$, i.e.\ a positive normalized functional on that algebra \cite{parthasarathy92}. In our paper this state will always be a vector state, i.e.\ of the form $\Phi(X)=\langle\Omega,X\Omega\rangle$ for some unit vector $\Omega\in H$. (Note that our inner products are linear in the second argument.)

If $X$ is a quantum random variable, i.e.\ an operator on a Hilbert space $H$, for which a functional calculus $C_b(\mathbb{C})\ni h\mapsto h(X)\in\mathcal{B}(H)$ for bounded continuous functions can be defined, then we call a probability measure $\mu$ on $\mathbb{C}$ the {\em distribution} of $X$ with respect to the vector state given by $\Omega\in H$, if
\[
\langle\Omega,h(X)\Omega\rangle = \int_\mathbb{C}h(x){\rm d}\mu(x)
\]
for all $h\in C_b(\mathbb{C})$. In this case we shall also write $\mu=\mathcal{L}(X,\Omega)$.

A densely defined, closed operator $X$ on a Hilbert space $H$ with domain ${\rm Dom}\,X$ is called {\em normal}, if the domains of $XX^*$ and $X^*X$ coincide and we have $XX^*=X^*X$ on this common domain.  

By the spectral theorem, an operator $X$ on a separable Hilbert space $H$ is normal if an only if there exists a $\sigma$-finite measure space $(E,\mathcal{E},\mu)$ and an $\mathcal{E}$-measurable function $\phi$ s.t.\ $X$ is unitarily equivalent to the operator $M_\phi$ of multiplication by $\phi$ on $L^2(E,\mathcal{E},\nu)$, i.e., there exists a unitary operator $U:H\to L^2(E,\mathcal{E},\nu)$ s.t.\ $UX\xi = \phi U\xi$ for all $\xi\in{\rm Dom}\,X$, see, e.g., \cite{conway85,pedersen89}. It follows that normal operators are quantum random variables in the sense above, since a functional calculus for them can be defined by $h(X)=U^*M_{h\circ\phi}U$. The distribution of $X$ w.r.t.\ to a unit vector $\Omega\in H$ is given by
\[
\mathcal{L}(X,\Omega) = \phi^{-1}\left(|U\Omega|^2\nu\right),
\]
since
\[
\langle \Omega, h(X)\Omega\rangle = \int_E h\circ\phi\, |U\Omega|^2{\rm d}\nu
\]
for all $h\in C_b(\mathbb{C})$.

To deal with unitary equivalence of possibly unbounded normal operators, we will use the following lemma.
\begin{lemma}\label{lem-unitary-eq}
Let $X$ and $X'$ be possibly unbounded normal operators on Hilbert spaces $H$ and $H'$. Assume that there exists a unitary operator $U:H\to H'$ such that
\begin{equation}\label{eq-unitary-lem-cond}
Uf(X) = f(X')U
\end{equation}
for any bounded continuous function $f$ on $\mathbb{C}$.

Then $UX = X'U$.
\end{lemma}
To prove this lemma, one can first prove the analogous statement for self-adjoint operators, e.g., using resolvents, and then apply the properties of the decomposition $X=A+iB$ of normal perators as a linear combination of two commuting self-adjoint operators, cf.\ \cite[Proposition 5.1.10]{pedersen89}. It is actually sufficient to require Condition \eqref{eq-unitary-lem-cond} for a much smaller class of functions, e.g.\ compactly supported real-valued $C^\infty$-functions.

\subsection{Nevanlinna theory and Cauchy transforms}\label{nevanlinna}

Denote by $\mathbb{C}^+=\{z\in\mathbb{C};{\rm Im}\,z>0\}$ and $\mathbb{C}^-=\{z\in\mathbb{C};{\rm Im}\,z<0\}$ the upper and lower half plane. For $\mu$ a probability measure on $\mathbb{R}$ and $z\in\mathbb{C}^+$, we define its {\em Cauchy transform} $G_\mu$ by
\[
G_\mu(z) = \int_{\mathbb{R}} \frac{1}{z-x}{\rm d}\mu(x)
\]
and its {\em reciprocal Cauchy transform} $F_\mu$ by
\[
F_\mu(z) = \frac{1}{G_\mu(z)}.
\]
Denote by $\mathcal{F}$ the following class of holomorphic
self-maps,
\[
\mathcal{F}=\left\{F:\mathbb{C}^+\to\mathbb{C}^+; F\mbox{
holomorphic and }\inf_{z\in\mathbb{C}^+} \frac{{\rm Im}\, F(z)}{{\rm
Im}\,z} = 1\right\}
\]
The map $\mu\mapsto F_\mu$ defines a bijection between the class $\mathcal{M}_1(\mathbb{R})$ of probability measures on $\mathbb{R}$ and $\mathcal{F}$, as follows from the following theorem.

\begin{theorem}\cite{maassen92}
Let $F: \mathbb{C}^+\to\mathbb{C}^+$ be holomorphic, then the
following are equivalent.
\begin{description}
\item[(i)] $\inf_{z\in\mathbb{C}^+} \frac{{\rm Im}\, F(z)}{{\rm
Im}\,z} = 1$;
\item[(ii)]
there exists a $\mu\in\mathcal{M}_1(\mathbb{R})$ such that
$F=F_\mu$.
\end{description}
Furthermore, $\mu$ is uniquely determined by $F$.
\end{theorem}

Similarly, for $\mu$ a probability measure on the unit circle
$\mathbb{T}=\{z\in\mathbb{C};|z|=1\}$ or on the positive half-line
$\mathbb{R}_+=\{x\in\mathbb{R};x\ge 0\}$, we define
\[
\psi_\mu(z) = \int \frac{xz}{1-xz}{\rm d}\mu
\]
and
\[
K_\mu(z)=\frac{\psi_\mu(z)}{1+\psi_\mu(z)}
\]
for $z\in\mathbb{C}\backslash{\rm supp}\,\mu$.

The map $\mu\mapsto K_\mu$ defines bijections between the class $\mathcal{M}_1(\mathbb{T})$ of probability measures on $\mathbb{T}$ and the class
\[
\mathcal{S}=\{K:\mathbb{D}\to\mathbb{D};K\mbox{ holomorphic and
}K(0)=0\},
\]
where $\mathbb{D}=\{z\in\mathbb{C};|z|<1\}$, and between the class $\mathcal{M}_1(\mathbb{R}_+)$ of probability measures on $\mathbb{R}_+$ and the class
\[
\mathcal{P}=\left\{K:\mathbb{C}\backslash\mathbb{R}_+\to
\mathbb{C}\backslash\mathbb{R}_+; \begin{array}{c} K \mbox{
holomorphic},\lim_{t\nearrow 0} K(t)=0,
K(\overline{z})=\overline{K(z)}, \\ \pi\ge\arg K(z)\ge\arg z\mbox{
for all }z\in \mathbb{C}^+\end{array}\right\},
\]
cf.\ \cite{belinschi+bercovici05} and the references therein.

In the following, if $X$ is an operator with distribution $\mu=\mathcal{L}(X,\Omega)$ w.r.t.\ $\Omega$, then we will write $G_X$, $F_X$, $\Psi_X$ or $K_X$ instead of $G_{\mathcal{L}(X,\Omega)}$, $F_{\mathcal{L}(X,\Omega)}$, $\psi_{\mathcal{L}(X,\Omega)}$, or $K_{\mathcal{L}(X,\Omega)}$ for the transforms of the distribution of $X$.

\subsection{Free convolutions}\label{free-conv}

By $\mathbb{A}_k$ we call denote the set of alternating $k$-tuples of 1's and 2's, i.e.
\[
\mathbb{A}_k= \big\{ (\varepsilon_1,\ldots,\varepsilon_k)\in \{1,2\}^k; \varepsilon_1\not=\varepsilon_2\not=\ldots\not=\varepsilon_k\big\}.
\]

\begin{definition}
\cite{voiculescu86}
Let $\mathcal{A}_1,\mathcal{A}_2\subseteq\mathcal{B}(H)$ be two $*$-algebras of bounded operators on a Hilbert space and assume $\mathbf{1}\in\mathcal{A}_i$, $i=1,2$. Let $\Omega$ be a unit vector in $H$ and denote by $\Phi$ the vector state associated to $\Omega$. We say that $\mathcal{A}_1$ and $\mathcal{A}_2$ are {\em free}, if we have
\[
\Phi(X_1\cdots X_k) =0
\]
for all $k\ge 1$, $\varepsilon \in\mathbb{A}_k$, $X_1\in\mathcal{A}_{\varepsilon_1},\ldots,X_k\in\mathcal{A}_{\varepsilon_k}$ such that
\[
\Phi(X_1)=\cdots=\Phi(X_k)=0.
\]
\end{definition}

Two normal operators $X$ and $Y$ are called {\em free}, if the algebras ${\rm alg}(X)=\{h(X);h\in C_b(\mathbb{C})\}$ and ${\rm alg}(Y)=\{h(Y);h\in C_b(\mathbb{C})\}$ they generate are free.

\begin{theorem}\label{thm-parallelogramm}
\cite{maassen92,chistyakov+goetze05,chistyakov+goetze06}
Let $\mu$ and $\nu$ be two probability measures on the real line, with reciprocal Cauchy transforms $F_\mu$ and $F_\nu$. Then there exist unique functions $Z_1,Z_2\in\mathcal{F}$ such that
\[
F_\mu\big(Z_1(z)\big) = F_\nu\big(Z_2(z)\big)= Z_1(z)+Z_2(z)-z
\]
for all $z\in\mathbb{C}^+$.
\end{theorem}
The function $F=F_\mu\circ Z_1 = F_\nu\circ Z_2$ also belongs to $\mathcal{F}$ and is therefore the the reciprocal Cauchy transform of some probability measure $\lambda$. One defines the additive free convolution of $\mu$ and $\nu$ as this unique probability measure and writes $\mu\boxplus\nu=\lambda$. This is justified by the following theorem.

\begin{theorem}
\cite{maassen92,bercovici+voiculescu93}
Let $X$ and $Y$ be two self-adjoint operators on some Hilbert space $H$ that are free w.r.t.\ some unit vector $\Omega\in H$. If $\Omega$ is cyclic, i.e.\ if
\[
\overline{{\rm alg}\{h(X),h(Y);h\in C_b(\mathbb{R})\}\Omega}=H.
\]
then $X+Y$ is essentially self-adjoint and the distribution w.r.t.\ $\Omega$ of its closure is equal to the additive free convolution of the distributions of $X$ and $Y$ w.r.t.\ to $\Omega$, i.e.
\[
\mathcal{L}(X+Y,\Omega) = \mathcal{L}(X,\Omega)\boxplus\mathcal{L}(Y,\Omega).
\]
\end{theorem}

There exist analogous results for the multiplicative convolutions of probability measures on the unit circle and the positive half-line, cf.\ \cite{maassen92,bercovici+voiculescu93,chistyakov+goetze05,chistyakov+goetze06}

\begin{theorem}
\item[(i)]
Let $\mu$ and $\nu$ be two probability measures  on the unit circle with transforms $K_\mu$ and $K_\nu$ and whose first moments do not vanish, $\int_\mathbb{T}x{\rm d}\mu(x)\not=0$, $\int_\mathbb{T}x{\rm d}\nu(x)\not=0$. Then there exist unique functions $Z_1,Z_2\in\mathcal{S}$
such that
\[
K_\mu\big(Z_1(z)\big) = K_\nu\big(Z_2(z)\big)=\frac{Z_1(z)Z_2(z)}{z}
\]
for all $z\in\mathbb{D}\backslash\{0\}$. The multiplicative free convolution $\lambda=\mu\boxtimes\nu$ is defined as the unique probability measure $\lambda$ with transform $K_\lambda=K_\mu\circ Z_1 = K_\nu\circ Z_2$.
\item[(ii)]
Let $U$ and $V$ be two unitary operators on some Hilbert space $H$ that are free w.r.t.\ some unit vector $\Omega\in H$. Then the products $UV$ and $VU$ are also unitary and their distributions w.r.t.\ to $\Omega$ are equal to the free convolution of the distributions of $U$ and $V$ w.r.t.\ $\Omega$, i.e.\, i.e.\
\[
\mathcal{L}(UV,\Omega) = \mathcal{L}(VU,\Omega) = \mathcal{L}(U,\Omega)\boxtimes\mathcal{L}(V,\Omega).
\]
\end{theorem}

\begin{theorem}
\item[(i)]
Let $\mu$ and $\nu$ be two probability measures  on the positive half-line such that $\mu\not=\delta_0$, $\nu\not=\delta_0$ and denote their transforms by $K_\mu$ and $K_\nu$. Then there exist unique functions $Z_1,Z_2\in\mathcal{P}$ such that
\[
K_\mu\big(Z_1(z)\big) = K_\nu\big(Z_2(z)\big)=\frac{Z_1(z)Z_2(z)}{z}
\]
for all $z\in\mathbb{C}\backslash\mathbb{R}_+$. The multiplicative free convolution $\lambda=\mu\boxtimes\nu$ is defined as the unique probability measure $\lambda$ with transform $K_\lambda=K_\mu\circ Z_1 = K_\nu\circ Z_2$.
\item[(ii)]
Let $X$ and $Y$ be two positive operators on some Hilbert space $H$ that are free w.r.t.\ some unit vector $\Omega\in H$.  Assume furthermore that $\Omega$ is cyclic, i.e. that
\[
\overline{{\rm alg}\{h(X),h(Y);h\in C_b(\mathbb{R})\}\Omega}=H.
\]
Then the products $\sqrt{X}Y\sqrt{X}$ and $\sqrt{Y}X\sqrt{Y}$ are essentially self-adjoint and positive, and their distributions w.r.t.\ to $\Omega$ are equal to the free convolution of the distributions of $X$ and $Y$ w.r.t.\ $\Omega$, i.e.\, i.e.\
\[
\mathcal{L}(\sqrt{X}Y\sqrt{X},\Omega) = \mathcal{L}(\sqrt{Y}X\sqrt{Y},\Omega) = \mathcal{L}(X,\Omega)\boxtimes\mathcal{L}(Y,\Omega).
\]
\end{theorem}

\section{Monotone Convolutions}

\begin{definition}\cite{muraki00}
Let $\mathcal{A}_1,\mathcal{A}_2\subset\mathcal{B}(H)$ be two
$*$-algebras of bounded operators on a Hilbert space $H$, and let
$\Omega\in H$ be a unit vector. We say that $\mathcal{A}_1$ and
$\mathcal{A}_2$ are {\em monotonically independent} w.r.t.\
$\Omega$, if we have
\[
\langle \Omega, X_1X_2\cdots X_k\Omega\rangle = \left\langle\Omega,\prod
_{\kappa:\varepsilon_\kappa=1} X_\kappa\Omega\right\rangle\prod
_{\kappa:\varepsilon_\kappa=2} \langle \Omega, X_\kappa\Omega\rangle
\]
for all $k\in\mathbb{N}$, $\varepsilon\in\mathbb{A}_k$,
$X_1\in\mathcal{A}_{\varepsilon_1},\ldots,X_k\in\mathcal{A}_{\varepsilon_k}$.
\end{definition}

\begin{remark}\label{rem-mon-cyclic}
\item[(a)]
Note that this notion depends on the order, i.e.\ if $\mathcal{A}_1$ and $\mathcal{A}_2$ are monotonically independent, then this does {\em not} imply that $\mathcal{A}_2$ and $\mathcal{A}_1$ are monotonically independent. In fact, if  $\mathcal{A}_1$ and $\mathcal{A}_2$ are monotonically independent and $\mathcal{A}_2$ and $\mathcal{A}_1$ are also monotonically independent, and $\Phi(\cdot)=\langle\Omega,\,\cdot\,\Omega\rangle$ does not vanish on one of the algebras, then restrictions of $\Phi$ to $\mathcal{A}_1$ and $\mathcal{A}_2$ have to be homomorphisms. To prove this for the restriction to, e.g., $\mathcal{A}_1$, take an element $Y\in\mathcal{A}_2$ such that $\Phi(Y)\not=0$, then
\[
\Phi(X_1X_2) = \frac{\Phi(X_1YX_2)}{\Phi(Y)} = \Phi(X_1)\Phi(X_2)
\]
for all $X_1,X_2\in\mathcal{A}_1$.
\item[(b)]
The algebras are not required to be unital. If  $\mathcal{A}_1$ is
unital, then the restriction of
$\Phi(\cdot)=\langle\Omega,\,\cdot\,\Omega\rangle$ to
$\mathcal{A}_2$ has to be a homomorphism, since monotone
independence implies
\[
\langle \Omega, XY\Omega\rangle =  \langle \Omega,
X\mathbf{1}Y\Omega\rangle  = \langle \Omega, X\Omega\rangle\langle
\Omega,Y\Omega\rangle
\]
for $X,Y\in\mathcal{A}_2$.
\item[(c)]
In the definition of monotone independence the condition
\[
XYZ= \langle\Omega,Y\Omega\rangle XZ
\]
for all $X,Z\in\mathcal{A}_1$, $Y\in\mathcal{A}_2$ is often also
imposed. If the state vector $\Omega$ is cyclic for the algebra
generated by $\mathcal{A}_1$ and $\mathcal{A}_2$, then this is
automatically satisfied. Let $X_1,X_3,\ldots,Z_1,Z_3,\ldots\in\mathcal{A}_1$ and $Y,X_2,X_4,\ldots,Z_2,Z_4,\ldots\in\mathcal{A}_2$, then
\begin{eqnarray*}
&& \langle X_1\cdots X_n\Omega, YZ_1\cdots Z_m\Omega\rangle = \langle \Omega,X_n^*\cdots X_1^* YZ_1\cdots Z_m\Omega\rangle \\
&=& \langle \Omega, Y\Omega\rangle \prod_{k \mbox{ even}}\langle \Omega, X^*_k\Omega\rangle\prod_{\ell \mbox{ even}}\langle \Omega, Z_\ell\Omega\rangle \langle X_1X_3\cdots\Omega,Z_1Z_3\cdots\Omega \\
&=&\langle \Omega, Y\Omega\rangle\langle X_1\cdots X_n\Omega, Z_1\cdots Z_m\Omega\rangle,
\end{eqnarray*}
for all $n,m\ge 1$, i.e., $X_1^*YZ_1$ and $\langle\Omega, Y\Omega\rangle X_1^*Z_1$ coincide on the subspace generated by $\mathcal{A}_1$ and $\mathcal{A}_2$ from $\Omega$.
\end{remark}

\begin{definition}\label{def-mon-indep}
Let $X$ and $Y$ be two normal operators on a Hilbert space $H$,
not necessarily bounded. We say that $X$ and $Y$ are {\em
monotonically independent} w.r.t.\ $\Omega$, if the $*$-algebras
${\rm alg}_0(X)=\{h(X);h\in C_b( \mathbb{C}),h(0)=0\}$ and ${\rm
alg}_0(Y)=\{h(Y);h\in C_b( \mathbb{C}),h(0)=0\}$ are monotonically
independent w.r.t.\ $\Omega$.
\end{definition}

Let us now introduce the model we shall use for calculations with monotonically independent operators.

\begin{proposition}\label{prop-mon-model}
Let $\mu,\nu$ be two probability measures on $\mathbb{C}$ and define normal operators $X$ and $Y$ on $L^2(\mathbb{C}\times\mathbb{C},\mu\otimes\nu)$ by
\begin{eqnarray*}
{\rm Dom}\, X &=& \left\{\psi\in L^2(\mathbb{C}\times\mathbb{C},\mu\otimes\nu); \int_{\mathbb{C}}\left|x\int_{\mathbb{C}}\psi(x,y){\rm d}\nu(y)\right|^2{\rm d}\mu(x)<\infty\right\}, \\
{\rm Dom}\, Y &=& \left\{\psi\in L^2(\mathbb{C}\times\mathbb{C},\mu\otimes\nu); \int_{\mathbb{C}\times\mathbb{C}}|y\psi(x,y)|^2{\rm d}\mu\otimes\nu(x,y)<\infty\right\},
\end{eqnarray*}
\begin{eqnarray*}
(X\psi)(x,y) &=& x\int_{\mathbb{C}}\psi(x,y'){\rm d}\nu(y'), \\
(X\psi)(x,y) &=& y\psi(x,y).
\end{eqnarray*}
Then $\mathcal{L}(X,\mathbf{1}) = \mu$, $\mathcal{L}(Y,\mathbf{1}) = \nu$, and $X$ and $Y$ are monotonically independent w.r.t.\ the constant function $\mathbf{1}$.
\end{proposition}
\begin{proof}
Denote by $P_2$ the orthogonal projection onto the space of functions in $L^2(\mathbb{C}\times\mathbb{C},\mu\otimes\nu)$ which do not depend on the second variable, and by $M_x$ multiplication by the first variable, then $X=M_xP_2$. This operator is normal, we have
\[
h(X)\psi(x,y) = \big(h(x)-h(0)\big)\int_{\mathbb{C}}\psi(x,y){\rm d}\nu(y) + h(0)\psi(x,y)
\]
and  $\langle \mathbf{1},h(X)\mathbf{1}\rangle = \int_{\mathbb{C}} h(x){\rm d}\mu(x)$ for all $h\in C_b(\mathbb{C})$, i.e.\ $\mathcal{L}(X,\mathbf{1})=\mu$. The operator $Y$ is multiplication by the second variable, it is clearly normal. We have
\[
h(Y)\psi(x,y)= h(y)\psi(x,y)
\]
and $\langle \mathbf{1},h(Y)\mathbf{1}\rangle = \int_{\mathbb{C}} h(y){\rm d}\nu(y)$ for all $h\in C_b(\mathbb{C})$, i.e.\ $\mathcal{L}(Y,\mathbf{1})=\nu$.

Let $f_1,\ldots,f_n,g_1,\ldots,g_n\in C_b(\mathbb{C})$, $f_1(0)=\cdots=f_n(0)=0$. Then
\[
f_n(X)g_{n-1}(Y)\cdots g_1(Y)f_1(X)\mathbf{1}=\prod_{k=1}^{n-1}\int_{\mathbb{C}}g_k(y){\rm d}\nu(y)\, f_1\cdots f_n
\]
and
\begin{eqnarray*}
\langle\mathbf{1},f_n(X)g_{n-1}(Y)\cdots g_1(Y)f_1(X)\mathbf{1}\rangle &=& \prod_{k=1}\int_{\mathbb{C}}g_k(y){\rm d}\nu(y)\int_{\mathbb{C}} f_1(x)\cdots f_n(x){\rm d}\mu(x) \\
&=& \prod_{k=1}^{n-1}\langle\mathbf{1},g_k(Y)\mathbf{1}\rangle\langle\mathbf{1} f_1(X)\cdots f_n(X)\mathbf{1}\rangle,
\end{eqnarray*}
i.e.\ the condition for monotone independence is satisfied in this case.
Similarly one checks the expectation of $g_n(Y)f_n(X)\cdots g_1(Y)f_1(X)$, $f_n(X)g_n(Y)\cdots f_1(X)g_1(Y)$, and $g_n(Y)f_{n-1}(X)\cdots f_1(X)g_1(Y)$.
\end{proof}

The following theorem shows that any pair of monotonically independent normal operators can be reduced to this model.

\begin{theorem}\label{thm-mon-model}
Let $X$ and $Y$ be two normal operators on a Hilbert space
$H$ that are monotonically independent with respect to $\Omega\in H$
and let $\mu=\mathcal{L}(X,\Omega)$, $\nu=\mathcal{L}(Y,\Omega)$.

Then there exists an isometry $W:L^2(\mathbb{C}\times\mathbb{C},\mu\otimes \nu)\to H$ such that
\begin{eqnarray}\label{eq-mon-model}
W^*h(X)W\psi(x,y) &=& \big(h(x)-h(0)\big)\int \psi(x,y){\rm d}\nu(y)
+
h(0)\psi(x,y), \\
W^*h(Y)W\psi(x,y) &=& h(y)\psi(x,y) \nonumber
\end{eqnarray}
for $x,y\in\mathbb{C}$, $\psi\in L^2(\mathbb{C}\times\mathbb{C},\mu\otimes \nu)\cong L^2(\sigma_X,\mu)\otimes L^2(\sigma_Y,\nu)$ and $h\in C_b(\mathbb{C})$.

We have $WL^2(\mathbb{C}\times\mathbb{C},\mu\otimes\nu) = \overline{{\rm alg}\{h(X),h(Y);h\in C_b( \mathbb{C})\}\Omega}$.

If the vector $\Omega\in H$ is cyclic for the algebra ${\rm
alg}(X,Y)={\rm alg}\{h(X),h(Y);h\in C_b( \mathbb{C})\}$
generated by $X$ and $Y$, then $W$ is unitary.
\end{theorem}
\begin{proof}
Define $W$ on simple tensors of bounded continuous functions by
\[
W f\otimes g = g(Y)f(X)\Omega
\]
for $f,g\in C_b(\mathbb{C})$. It follows from the monotone
independence of $X$ and $Y$ that this defines an isomorphism, since
\begin{eqnarray*}
\langle W f_1\otimes g_1, W f_2\otimes g_2\rangle &=& \langle
\Omega, f_1(X)^* g_1(Y)^* g_2(Y)f_2(X)\Omega\rangle \\
&=& \langle \Omega, f_1(X)^*f_2(X)\Omega\rangle \langle \Omega,
g_1(Y)^*g_2(Y)\Omega\rangle \\
&=& \int \overline{f_1(t)}f_2(t){\rm d}\mu(t) \int
\overline{g_1(t)}g_2(t){\rm d}\nu(t).
\end{eqnarray*}
Since $C_b(\mathbb{C})\otimes C_b(\mathbb{C})$ is dense in
$L^2(\mathbb{C}\times\mathbb{C},\mu\otimes \nu)$, $W$ extends to a
unique isomorphism on $L^2(\mathbb{C}\times\mathbb{C},\mu\otimes \nu)$.

The relations
\begin{multline*}
\langle W f_1\otimes g_1, h(X) W f_2\otimes g_2\rangle = \langle
\Omega, f_1(X)^* g_1(Y)^* h(X) g_2(Y)f_2(X)\Omega\rangle \\
=  \langle \Omega, f_1(X)^* \big(h(X)-h(0)\big)f_2(X)\Omega\rangle
\langle \Omega, g_1(Y)^*\Omega\rangle  \langle \Omega,
g_2(Y)\Omega\rangle \\
+h(0) \langle \Omega, f_1(X)^* g_1(Y)^*g_2(Y)f_2(X)\Omega\rangle
\\
=\langle\Omega,g_2(Y)\Omega\rangle\left\langle W f_1\otimes g_1, W
\big((h-h(0)1\big)f_1\otimes
 1\right\rangle +  h(0)\langle W f_1\otimes g_1, Wf_1\otimes g_2\rangle \\
= \left\langle W f_1\otimes g_1, W\left( \int g_2(y){\rm
d}\nu(y)(h-h(0)1)f_1\otimes 1 + h(0) f_2\otimes g_2\right)
\right\rangle
\end{multline*}
and
\begin{eqnarray*}
\langle W f_1\otimes g_1, h(Y) W f_2\otimes g_2\rangle &=& \langle
\Omega, f_1(X)^* g_1(Y)^* h(Y) g_2(Y)f_2(X)\Omega\rangle \\
&=&  \langle W f_1\otimes g_1, W f_2\otimes (h g_2)\rangle
\end{eqnarray*}
shows that we have the desired formulas for simple tensors of functions $f_1,f_2,g_1,g_2\in C_b(\mathbb{C})$. The general case follows by linearity and continuity. Remark \ref{rem-mon-cyclic}(c) implies
\begin{eqnarray*}
WL^2(\mathbb{C}\times\mathbb{C},\mu\otimes\nu) &=& \overline{{\rm span}\,\{g(Y)f(X)\Omega;f,g\in C_b(\mathbb{C})\}} \\
&=& \overline{{\rm alg}\{h(X),h(Y);h\in C_b( \mathbb{C})\}\Omega}.
\end{eqnarray*}
If $\Omega$ is cyclic, then $W$ is surjective and therefore unitary.
\end{proof}

\begin{remark}\label{rem-mon-reduction}
It follows that the joint law of two monotonically independent, normal operators is uniquely determined by their marginal distributions, in the sense that the restriction of $\Phi(\cdot)=\langle\Omega,\cdot\,\Omega\rangle$ to ${\rm alg}(X,Y)={\rm alg}\{h(X),h(Y);h\in C_b( \mathbb{C})\}$ is uniquely determined by $\mathcal{L}(X,\Omega)$ and $\mathcal{L}(Y,\Omega)$. But by Lemma \ref{lem-unitary-eq}, also computations for unbounded functions of $X$ and $Y$, e.g., concerning the operators $X+Y$ for self-adjoint $X$ and $Y$, or $\sqrt{X}Y\sqrt{Y}$ for positive $X$ and $Y$, reduce to the model introduced in Proposition \ref{prop-mon-model}.
\end{remark}

\subsection{Additive monotone convolution on
$\mathcal{M}_1(\mathbb{R})$}

\begin{definition} \cite{muraki00}
Let $\mu$ and $\nu$ be two probability measures on $\mathbb{R}$ with
reciprocal Cauchy transforms $F_\mu$ and $F_\nu$. Then we define the
additive monotone convolution $\lambda=\mu\triangleright\nu$ of
$\mu$ and $\nu$ as the unique probability measure on $\mathbb{R}$
with reciprocal Cauchy transform $F_\lambda = F_\mu \circ F_\nu$.
\end{definition}

It follows from Subsection \ref{nevanlinna} that the additive
monotone convolution is well-defined. Let us first recall some basic
properties of the additive monotone convolution.

\begin{proposition}\cite{muraki00}\label{prop-mon-add-properties}
The additive monotone convolution is associative and $*$-weakly continuous in both arguments. It is affine in the first argument and convolution from the right by a Dirac measure corresponds to translation, i.e. $\mu\triangleright\delta_x = T_x^{-1}\mu$ for $x\in\mathbb{R}$, where $T_x:\mathbb{R}\to\mathbb{R}$ is defined by $T_x(t)=t+x$.
\end{proposition}
This convolution is not commutative, i.e.\ in general we have $\mu\triangleright\nu\not=\nu\triangleright\mu$.

Let $x\in\mathbb{R}$ and $0\le p \le 1$. Then one can compute, e.g.,
\[
\delta_x\triangleright\big(p\delta_{1} + (1-p)\delta_{-1}\big) = q\delta_{z_1} + (1-q)\delta_{z_2}
\]
where
\begin{eqnarray*}
z_1 &=& \frac{1}{2}\left(x +\sqrt{x^2+4(2p-1)x+4}\right) , \\
z_2 &=& \frac{1}{2}\left(x -\sqrt{x^2+4(2p-1)x+4}\right) , \\
q &=& \frac{x+4p-2+\sqrt{x^2+4(2p-1)x+4}}{2\sqrt{x^2+4(2p-1)x+4}}.
\end{eqnarray*}
This example shows that convolution from the left by a Dirac mass is in general not equal to a translation and that the additive monotone convolution is not affine in the second argument.

Note that the continuity and the fact that the monotone convolution is affine in the first argument imply the following formula
\begin{equation}\label{eq-int-formula}
\mu\triangleright\nu = \int_\mathbb{R} \delta_x\triangleright \nu\, {\rm d}\mu(x)
\end{equation}
for all $\mu,\nu\in\mathcal{M}_1(\mathbb{R})$.

The following proposition is the key to treating the additive monotone convolution for general probability measures on $\mathbb{R}$.

\begin{proposition}\label{prop-mon-add-resolvent}
Let $\mu$ and $\nu$ be two probability measures on $\mathbb{R}$ and denote by $M_x$ and $M_y$ the self-adjoint operators on $L^2(\mathbb{R}\times\mathbb{R},\mu\otimes\nu)$ defined by multiplication with the coordinate functions. Denote by $P_2$ the orthogonal projection onto the subspace of functions which do not depend on the second coordinate, $L^2(\mathbb{R}\times\mathbb{R},\mu\otimes\nu)\ni\psi\mapsto \int_\mathbb{R}\psi(\cdot,y){\rm d}\nu(y)\in L^2(\mathbb{R}\times\mathbb{R},\mu\otimes\nu)$. Then $M_xP_2=P_2M_x$ and $M_y$ are self-adjoint and monotonically independent w.r.t.\ the constant function and the operator $z-M_xP_2-M_y$ has a bounded inverse for all $z\in\mathbb{C}\backslash\mathbb{R}$, given by
\begin{equation}\label{eq-mon-add-resolvent}
\left((z-M_xP_2-M_y)^{-1}\psi\right)(x,y) = \frac{\psi(x,y)}{z-y} + \frac{x \int_\mathbb{R} \frac{\psi(x,y')}{z-y'}{\rm d}\nu(y')}{(z-y)(1-xG_\nu(z))}.
\end{equation}
\end{proposition}
\begin{proof}
$M_xP_2$ and $M_y$ are monotonically independent by \ref{prop-mon-model}.

The first term on the right-hand-side of Equation \eqref{eq-mon-add-resolvent} is obtained from $\psi$ by multiplication with a bounded function, the second by composition of multiplications with bounded functions and the projection $P_2$. Equation \eqref{eq-mon-add-resolvent} therefore clearly defines a bounded operator. To check that it is indeed the inverse of $z-M_xP_2-M_y$ is straightforward,
\begin{multline*}
(z-M_xP_2-M_y) \left(\frac{\psi(x,y)}{z-y} + \frac{x \int_\mathbb{R} \frac{\psi(x,y')}{z-y'}{\rm d}\nu(y')}{(z-y)\big(1-xG_\nu(z)\big)}\right) \\
= \psi(x,y) + \frac{x\int\frac{\psi(x,y')}{z-y'}{\rm d}\nu(y')}{1-xG_\nu(z)} - x\int_\mathbb{R} \frac{\psi(x,y')}{z-y'}{\rm d}\nu(y')-x\int_\mathbb{R}\frac{x \int_\mathbb{R} \frac{\psi(x,y'')}{z-y''}{\rm d}\nu(y'')}{(z-y')\big(1-xG_\nu(z)\big)}{\rm d}\nu(y') \\
= \psi(x,y) +\frac{\Big((z-y)-(z-y)\big(1-xG_\nu(z)\big) - xG_\nu(z)(z-y)\Big)x\int_\mathbb{R}\frac{\psi(x,y')}{z-y'}{\rm d}\nu(y')}{(z-y)\big(1-xG_\nu(z)\big)} \\
 = \psi(x,y)
\end{multline*}
\end{proof}

\begin{theorem}\label{thm-mon-add-conv}
Let $X$ and $Y$ be two self-adjoint operators on a Hilbert space $H$ that are monotonically independent w.r.t.\ to a unit vector $\Omega\in H$. Assume furthermore that $\Omega$ is cyclic, i.e. that
\[
\overline{{\rm alg}\{h(X),h(Y);h\in C_b(\mathbb{R})\}\Omega}=H.
\]
Then $X+Y$ is essentially self-adjoint and the distribution w.r.t.\ $\Omega$ of its closure is equal to the additive monotone convolution of the distributions of $X$ and $Y$ w.r.t.\ to $\Omega$, i.e.
\[
\mathcal{L}(X+Y,\Omega) = \mathcal{L}(X,\Omega)\triangleright\mathcal{L}(Y,\Omega).
\]
\end{theorem}
\begin{proof}
Let $\mu=\mathcal{L}(X,\Omega)$, $\nu=\mathcal{L}(Y,\Omega)$.

By Theorem \ref{thm-mon-model} and Lemma \ref{lem-unitary-eq} it is sufficient to consider the case where $X$ and $Y$ are given by Proposition \ref{prop-mon-model}. Proposition \ref{prop-mon-add-resolvent} shows that $z-X-Y$ admits a bounded inverse and therefore that ${\rm Ran}\,(z-X-Y)$ is dense for $z\in\mathbb{C}\backslash\mathbb{R}$. By \cite[Theorem VIII.3]{reed+simon80} this is equivalent to $X+Y$ being essentially self-adjoint.

Using Equation \eqref{eq-mon-add-resolvent}, we can compute the Cauchy transform of the distribution of the closure of $X+Y$. Let $z\in\mathbb{C}^+$, then we have
\begin{multline*}
G_{X+Y}(z) = \langle \Omega, (z-X-Y)^{-1}\Omega\rangle - = \left\langle \mathbf{1}, (z-M_xP_2-M_y)^{-1}\mathbf{1}\right\rangle \\
=\left\langle\mathbf{1},\frac{1}{z-y} + \frac{xG_\nu(z)}{(z-y)(1-xG_\nu(z))}\right\rangle = \int_{\mathbb{R}\times\mathbb{R}} \frac{1}{(z-y)(1-xG_\nu(z))}{\rm d}\mu\otimes\nu \\
= \int_\mathbb{R} \frac{G_\nu(z)}{1-xG_\nu(z)}{\rm d}\mu(x) = G_\mu\left(\frac{1}{G_\nu(z)}\right)=G_\mu\big(F_\nu(z)\big),
\end{multline*}
or
\[
F_{X+Y}(z) = \frac{1}{G_{X+Y}(z)} = \frac{1}{G_\mu\big(F_\nu(z)\big)} = F_\mu\big(F_\nu(z)\big)= F_{\mu\triangleright\nu}(z).
\]
\end{proof}

\subsection{Multiplicative monotone convolution on
$\mathcal{M}_1(\mathbb{R}_+)$}

\begin{definition}\cite{bercovici04}
Let $\mu$ and $\nu$ be two probability measures on the positive half-line $\mathbb{R}_+$ with transforms $K_\mu$ and $K_\nu$. Then the multiplicative monotone convolution of $\mu$ and $\nu$ is defined as the unique probability measure $\lambda=\mu\mmultg\nu$ on $\mathbb{R}_+$ with transform $K_\lambda=K_\mu\circ K_\nu$.
\end{definition}

It follows from Subsection \ref{nevanlinna} that the multiplicative monotone convolution on $\mathcal{M}_1(\mathbb{R}_+)$ is well-defined.

Let us first recall some basic properties of the multiplicative monotone convolution.

\begin{proposition}\label{prop-mon-mult1-properties}
The multiplicative monotone convolution $\mathcal{M}_1(\mathbb{R}_+)$ is associative and $*$-weakly continuous in both arguments. It is affine in the first argument and convolution from the right by a Dirac measure corresponds to dilation, i.e. $\mu\mmultg\delta_\alpha = D_\alpha^{-1}\mu$ for $\alpha\in\mathbb{R}_+$, where $D_\alpha:\mathbb{R}_+\to\mathbb{R}_+$ is defined by $D_\alpha(t)=\alpha t$.
\end{proposition}
This convolution is not commutative, i.e.\ in general we have $\mu\mmultg\nu\not=\nu\mmultg\mu$. As in the additive case is not affine in the second argument, either, and convolution from the left by a Dirac mass is in general not equal to a dilation.

We want to extend \cite[Corollary 4.3]{franz05a} to unbounded positive operators, i.e.\ we want to show that if $X$ and $Y$ are two positive operators such that $X-\mathbf{1}$ and $Y$ are monotonically independent, then the distribution of $\sqrt{X}Y\sqrt{X}$ is equal to the multiplicative monotone convolution of the distributions of $X$ and $Y$. By Theorem \ref{thm-mon-model}, it is sufficient to do the calculations for the case where $X$ and $Y$ are constructed from multiplication with the coordinate functions and the projection $P_2$.

\begin{proposition}\label{prop-mon-mult1-resolvent}
Let $\mu$ and $\nu$ be two probability measures on $\mathbb{R}_+$, $\nu\not=\delta_0$, and let $M_y$ be the self-adjoint operator on $L^2(\mathbb{R}_+\times\mathbb{R}_+,\mu\otimes\nu)$ defined by multiplication with the coordinate function $(x,y)\mapsto y$. Define $S_x$ on $L^2(\mathbb{R}_+\times\mathbb{R}_+,\mu\otimes\nu)$ by
\begin{eqnarray} \label{eq-sx-definition} \\
{\rm Dom}\,S_x &=& \left\{\psi\in L^2(\mathbb{R}_+\times\mathbb{R}_+,\mu\otimes\nu); \int_{\mathbb{R}_+} x\psi(x,y){\rm d}\nu(y)\in L^2(\mathbb{R}_+,\mu)\right\} , \nonumber \\
(S_x\psi)(x,y) &=& (x-1)\int_{\mathbb{R}_+} \psi(x,y){\rm d}\nu(y) + \psi(x,y) \nonumber
\end{eqnarray}
Then $S_x-\mathbf{1}$ and $M_y$ are monotonically independent w.r.t.\ to the constant function and the operator $z-\sqrt{S_x}M_y\sqrt{S_x}$ has a bounded inverse for all $z\in\mathbb{C}\backslash\mathbb{R}$, given by
\begin{equation}\label{eq-mon-mult1-resolvent}
\left((z-\sqrt{S_x}M_y\sqrt{S_x})^{-1}\psi\right)(x,y) =\frac{\psi(x,y)+g(x)}{z-y} + h(x).
\end{equation}
where
\begin{eqnarray*}
g(x) &=& \frac{\sqrt{x}-x}{(1-x)zG_\nu(z)+x}\int_{\mathbb{R}_+} \psi(x,y){\rm d}\nu(y) \\
&& + \frac{z(x-1)}{(1-x)zG_\nu(z)+x}\int_{\mathbb{R}_+}\frac{\psi(x,y)}{z-y}{\rm d}\nu(y), \\
h(x) &=& \frac{(\sqrt{x}-1)^2G_\nu(z)}{(1-x)zG_\nu(z)+x}\int_{\mathbb{R}_+}\psi(x,y){\rm d}\nu(y) \\
&& +\frac{\sqrt{x}-x}{(1-x)zG_\nu(z)+x}\int_{\mathbb{R}_+}\frac{\psi(x,y)}{z-y}{\rm d}\nu(y).
\end{eqnarray*}
\end{proposition}
\begin{proof}
Fix $z\in\mathbb{C}^+$. Let $x>0$, then
\[
{\rm Im}\frac{z}{z-x} = - \frac{x{\rm Im}\,z}{({\rm Re}\,z-x)^2+({\rm Im}\,z)^2}< 0,
\]
and therefore
\[
{\rm Im}\,zG_\nu(z) = {\rm Im}\int_{\mathbb{R}_+}\frac{z}{z-x}{\rm d}\nu(x) <0.
\]
Similarly, we get ${\rm Im}\,zG_\nu(z)>0$ for $z\in\mathbb{C}^-$. It follows that the functions in front of the integrals in the definitions of $g$ and $h$ are bounded as functions of $x$, and therefore $g$ and $h$ are square-integrable. Since $\frac{1}{z-y}$ is bounded, too, we see that Equation \eqref{eq-mon-mult1-resolvent} defines a bounded operator.

Let us now check that it is the inverse of $z-\sqrt{S_x}M_y\sqrt{S_x}$.

Using the notation of the previous subsection, we can write $S_x$ also as $S_x=M_{x-1}P_2+\mathbf{1}=M_xP_2 + P_2^\perp$, where $P_2^\perp$ is the projection onto the orthogonal complement of the subspace of functions which do not depend on $y$. Its square root can be written as $\sqrt{S_x}= M_{\sqrt{x}}P_2 + P_2^\perp = M_{\sqrt{x}-1}P_2 +\mathbf{1}$, it acts as
\[
\left(\sqrt{S_x}\psi\right)(x,y) =  \left(\sqrt{x}-1\right)\int_{\mathbb{R}_+} \psi(x,y){\rm d}\nu(y) + \psi(x,y)
\]
on a function $\psi\in{\rm Dom}\,\sqrt{S_x}\subseteq L^2(\mathbb{R}_+\times\mathbb{R}_+,\mu\otimes\nu)$.

Since $h$ does not depend on $y$, we have $\sqrt{S_x}h=\sqrt{x}h$. For $g$ we get
\begin{eqnarray*}
\left(\sqrt{S_x}\frac{g}{z-y}\right)(x) &=& (\sqrt{x}-1)\int_{\mathbb{R}_+}\frac{g(x)}{z-y}{\rm d}\nu(y) + \frac{g(x)}{z-y} \\
&=& \left((\sqrt{x}-1)G_\nu(z)+\frac{1}{z-y}\right) g(x).
\end{eqnarray*}
Set $\varphi=\frac{\psi+g}{z-y}+h$. Applying $\sqrt{S_x}$ to $\varphi$, we get
\begin{eqnarray*}
\left(\sqrt{S_x}\varphi\right)(x,y) &=& \frac{\psi(x,y)}{z-y}+\frac{\sqrt{x}-x}{(z-y)\big((1-x)zG_\nu(z)+x\big)}\int_{\mathbb{R}_+}\psi(x,y){\rm d}\nu(y) \\
&& +\frac{z(x-1)}{(z-y)\big((1-x)zG_\nu(z)+x\big)}\int_{\mathbb{R}_+}\frac{\psi(x,y)}{z-y}{\rm d}\nu(y) \\
&=& \frac{\psi(x,y)+g(x)}{z-y}.
\end{eqnarray*}
{}From this we get
\begin{eqnarray*}
\left(\left(z-\sqrt{S_x}M_y\sqrt{S_x}\right)\varphi\right)(x,y) &=& \psi(x,y)
\end{eqnarray*}
after some tedious, but straightforward computation.
\end{proof}

\begin{remark}
It $\nu=\delta_0$, then $M_y=0$ on $L^2(\mathbb{R}_+\times\mathbb{R}_+,\mu\otimes\nu)$, and therefore $\sqrt{S_x}M_y\sqrt{S_x}=0$. This is of course a positive operator, and its distribution is $\delta_0$.
\end{remark}

\begin{theorem}\label{thm-mon-mult}
Let $X$ and $Y$ be two positive self-adjoint operators on a Hilbert space $H$ such that $X-\mathbf{1}$ and $Y$ are monotonically independent w.r.t.\ to a unit vector $\Omega\in H$. Assume furthermore that $\Omega$ is cyclic, i.e.
\[
\overline{{\rm alg}\{h(X),h(Y);h\in C_b(\mathbb{R}_+)\}\Omega}=H.
\]
Then $\sqrt{X}Y\sqrt{X}$ is essentially self-adjoint and the distribution w.r.t.\ $\Omega$ of its closure is equal to the multiplicative monotone convolution of the distributions of $X$ and $Y$ w.r.t.\ $\Omega$, i.e.
\[
\mathcal{L}\left(\sqrt{X}Y\sqrt{X},\Omega\right)=\mathcal{L}(X,\Omega)\mmultg\mathcal{L}(Y,\Omega).
\]
\end{theorem}
\begin{proof}
Let $\mu=\mathcal{L}(X,\Omega)$, $\nu=\mathcal{L}(Y,\Omega)$.

By Theorem \ref{thm-mon-model} it is sufficient to consider the case $X=S_x$ and $Y=M_y$. In this case Proposition \ref{prop-mon-mult1-resolvent} shows that $z-\sqrt{X}Y\sqrt{X}$ has a bounded inverse for all $z\in\mathbb{C}\backslash\mathbb{R}$. This implies that ${\rm Ran}(z-\sqrt{X}Y\sqrt{X})$ is dense for all $z\in\mathbb{C}\backslash\mathbb{R}$ and that $\sqrt{X}Y\sqrt{X}$ is essentially self-adjoint, cf.\ \cite[Theorem VIII.3]{reed+simon80}.

Using Equation \eqref{eq-mon-mult1-resolvent}, we can compute the Cauchy transform of the distribution of the closure of $\sqrt{X}Y\sqrt{X}$. Let $z\in\mathbb{C}^+$, then we have
\begin{eqnarray*}
G_{\sqrt{X}Y\sqrt{X}}(z) &=& \left\langle\Omega,\left(z-\sqrt{X}Y\sqrt{X}\right)^{-1}\Omega\right\rangle = \left\langle \mathbf{1}, \left(z-\sqrt{S_x}M_y\sqrt{S_x}\right)^{-1}\mathbf{1}\right\rangle \\
&=& \left\langle \mathbf{1},\frac{1+g_1}{z-y}+h_1\right\rangle
\end{eqnarray*}
where
\begin{eqnarray*}
g_1(x) &=&\frac{\sqrt{x}-x+(x-1)zG_\nu(x)}{(1-x)zG_\nu(z)+x} = \frac{\sqrt{x}}{(1-x)zG_\nu(z)+x} -1, \\
h_1(x) &=&\frac{(1-\sqrt{x})G_\nu(z)}{(1-x)zG_\nu(z)+x}.
\end{eqnarray*}
Therefore
\begin{eqnarray}
G_{\sqrt{X}Y\sqrt{X}}(z) &=& \int_{\mathbb{R}_+\times\mathbb{R}_+} \left(\frac{1+g_1(x)}{z-y} + h_1(x)\right){\rm d}\mu\otimes\nu(x,y) \nonumber \\
&=& \int_{\mathbb{R}_+} \frac{G_\nu(z)}{(1-x)zG_\nu(z)+x}{\rm d}\mu(x) = \frac{G_\nu(z)}{zG_\nu(z)-1}G_\mu\left(\frac{zG_\nu(z)}{zG_\nu(z)-1}\right). \label{eq-xyx-yxy}
\end{eqnarray}
Using the relation
\[
G_\mu(z) = \frac{1}{z}\left(\psi_\mu\left(\frac{1}{z}\right)+1\right)
\]
to replace the Cauchy transforms by the $\psi$-transforms, this becomes
\[
\psi_{\sqrt{X}Y\sqrt{X}}\left(\frac{1}{z}\right) = \psi_\mu\left(\frac{\psi_\nu(1/z)}{\psi_\nu(1/z)+1}\right),
\]
or finally
\[
K_{\sqrt{X}Y\sqrt{X}}(z) = K_\mu\big(K_\nu(z)\big)=K_{\mu\mmultk\nu}(z).
\]
\end{proof}

\subsection{The ``other'' multiplicative convolution on $\mathcal{M}_1(\mathbb{R}_+)$}\label{subsec-mon-mult-other}

Let $X$ and $Y$ be two positive operators such that $X-\mathbf{1}$ and $Y$ are monotonically independent w.r.t.\ to some unit vector. We have just shown that the distribution of $\sqrt{X}Y\sqrt{X}$ is given by the multiplicative monotone convolution of the distributions of $X$ and $Y$, as defined by \cite{bercovici04}. But in \cite{franz05a} it is was already shown that in general this is {\em not} the case for $\sqrt{Y}X\sqrt{Y}$, which would be another obvious choice for constructing multiplicatively a positive operator out of $X$ and $Y$. It is possible to characterize the distribution of $\sqrt{Y}X\sqrt{Y}$ using the same methods as in the previous subsection. We will summarize the main results here, but omit the details of the calculations.

\begin{proposition}
Let $\mu$ and $\nu$ be two probability measures on $\mathbb{R}_+$, and $S_x$, $M_y$ as in Proposition \ref{prop-mon-mult1-resolvent}.

Then the operator $z-\sqrt{M_y}S_x\sqrt{M_y}$ has a bounded inverse for all $z\in\mathbb{C}\backslash\mathbb{R}$, given by
\begin{equation}
\left((z-\sqrt{M_y}S_x\sqrt{M_y})^{-1}\psi\right)(x,y) = \frac{\psi(x,y)+\sqrt{y}h(x)}{z-y},
\end{equation}
where
\[
h(x) = \frac{(x-1)\int_{\mathbb{R}_+} \frac{\sqrt{y}\psi(x,y)}{z-y}{\rm d}\nu(y)}{(1-x)zG_\nu(z)+x}.
\]
\end{proposition}

\begin{theorem}
Let $X$ and $Y$ be two positive self-adjoint operators on a Hilbert space $H$ such that $X-\mathbf{1}$ and $Y$ are monotonically independent w.r.t.\ to a unit vector $\Omega\in H$. Assume furthermore that $\Omega$ is cyclic, i.e.
\[
\overline{{\rm alg}\{h(X),h(Y);h\in C_b(\mathbb{R}_+)\}\Omega}=H.
\]
Then $\sqrt{Y}X\sqrt{Y}$ is positive and essentially self-adjoint. Denote by $\mu=\mathcal{L}(X,\Omega)$ and $\nu=\mathcal{L}(Y,\Omega)$ the distributions of $X$ and $Y$ w.r.t.\ $\Omega$, let $W_\nu(z)=G_{\sqrt{y}\nu}(z) = \int_{\mathbb{R}_+}\frac{\sqrt{y}}{z-y}{\rm d}\nu(y)$ for $y\in\mathbb{C}\backslash\mathbb{R}$. Then the distribution $\lambda=\mathcal{L}(\sqrt{Y}X\sqrt{Y},\Omega)$ of its closure w.r.t.\ $\Omega$ is characterized by its Cauchy transform
\begin{equation}\label{eq-mon-mult-alt}
G_\lambda(z) =G_\nu(z) -\frac{\big(W_\nu(z)\big)^2}{zG_\nu(z)-1} + \frac{\big(W_\nu(z)\big)^2}{\big(zG_\nu(z)-1\big)^2}G_\mu\left(\frac{zG_\nu(z)}{zG_\nu(z)-1}\right).
\end{equation}
\end{theorem}

\begin{remark}\label{rem-mon-mult-alt}
We can use Equation \eqref{eq-mon-mult-alt} to define an alternative multiplicative monotone convolution. Let $\mu$ and $\nu$ be two probability measures on $\mathbb{R}_+$, then we define $\lambda=\mu\mmultalt\nu$ as the unique probability measure $\lambda$ on $\mathbb{R}_+$ whose Cauchy transform is given by Equation \eqref{eq-mon-mult-alt}.

If $\nu=\delta_y$ is a Dirac mass, then $W_{\delta_y}=\frac{\sqrt{y}}{z-y}$, $\big(W_{\delta_y}(z)\big)^2=\frac{y}{(z-y)^2}=G_{\delta_y}(z)\big(zG_{\delta_y}(z)-1\big)$, and Equation \eqref{eq-mon-mult-alt} reduces to Equation \eqref{eq-xyx-yxy}. Therefore
\[
\mu\mmultalt\delta_y = \mu\mmultg\delta_y= D^{-1}_y\mu
\]
for all $y\in\mathbb{R}_+$ and $\mu\in\mathcal{M}_1(\mathbb{R}_+)$.

But in general the two convolutions are different, as was already stated in \cite{franz05a}. Actually, the convolution $\mmultalt$ is not even associative, as the following examples show. Let $x,y>0$, $0<p<1$, and set $X=\left(\begin{array}{cc} x & 0 \\ 0 & 1\end{array}\right)$, $Y=y \left(\begin{array}{cc} p & \sqrt{p(1-p)} \\ \sqrt{p(1-p)} & 1-p\end{array}\right)$, $\Omega=\left(\begin{array}{c} 1 \\ 0 \end{array}\right)$, then $X-\mathbf{1}$ and $Y$ are monotonically independent and $\mathcal{L}(X,\Omega)=\delta_x$, $\mathcal{L}(Y,\Omega=p\delta_0+(1-p)\delta_y$. Furthermore we get
\[
\delta_x\mmultalt\big(p\delta_0+(1-p)\delta_y) = \mathcal{L}\big(\sqrt{Y}X\sqrt{Y},\Omega\big) = p\delta_0 + (1-p)\delta_{y(xp+1-p)}.
\]
Therefore
\begin{eqnarray*}
(\delta_{x_1}\mmultalt\delta_{x_2})\mmultalt\big(p\delta_0+(1-p)\delta_y\big) &=&  \delta_{x_1x_2}\mmultalt\big(p\delta_0+(1-p)\delta_y\big) \\
&=& \big(p\delta_0+(1-p)\delta_{y(x_1x_2p+1-p)}\big) \\
&\not=& \delta_{x_1}\mmultalt\left(\delta_{x_2}\mmultalt\big(p\delta_0+(1-p)\delta_{y}\big)\right) \\
&=&  \delta_{x_1}\mmultalt\big(p\delta_0+(1-p)\delta_{y(x_2p+1-p)}\big)  \\
&=&\big(p\delta_0+(1-p)\delta_{y(x_1p+1-p)(x_2+1-p)}\big)  \\
\end{eqnarray*}
in general.
\end{remark}

\subsection{Multiplicative monotone convolution on $\mathcal{M}_1(\mathbb{T})$}\label{subsec-mult-mon-T}

\begin{definition}\cite{bercovici04}
Let $\mu$ and $\nu$ be two probability measure on the unit circle $\mathbb{T}$ with transforms $K_\mu$ and $K_\nu$. Then the multiplicative monotone convolution of $\mu$ and $\nu$ is defined as the unique probability measure $\lambda=\mu\mmultg\nu$ on $\mathbb{T}$ with transform $K_\lambda=K_\mu\circ K_\nu$.
\end{definition}

It follows from Subsection \ref{nevanlinna} that the multiplicative monotone convolution on $\mathcal{M}_1(\mathbb{T})$ is well-defined.

Let us first recall some basic properties of the multiplicative monotone convolution.

\begin{proposition}
The multiplicative monotone convolution on $\mathcal{M}_1(\mathbb{T})$ is associative and $*$-weakly continuous in both arguments. It is affine in the first argument and convolution from the right by a Dirac measure corresponds to rotation, i.e. $\mu\mmultg\delta_{e^{i\vartheta}} = R_\vartheta^{-1}\mu$ for $\vartheta\in[0,2\pi[$, where $R_\vartheta:\mathbb{T}\to\mathbb{T}$ is defined by $R_\vartheta(t)=e^{i\vartheta} t$.
\end{proposition}
This convolution is not commutative, i.e.\ in general we have $\mu\mmultg\nu\not=\nu\mmultg\mu$. As in the additive case is not affine in the second argument, either, and convolution from the left by a Dirac mass is in general not equal to a rotation.

Probability measures on the unit circle arise as distributions of unitary operators and they are completely characterized by their moments. Therefore the following theorem is a straightforward consequence of \cite{bercovici04} (see also \cite[Theorem 4.1 and Corollary 4.2]{franz05a}).

\begin{theorem}
Let $U$ and $V$ be two unitary operators on a Hilbert space $H$, $\Omega\in H$ a unit vector and assume furthermore that $U-\mathbf{1}$ and $V$ are monotonically independent w.r.t.\ $\Omega$. Then the products $UV$ and $VU$ are also unitary and their distribution w.r.t. $\Omega$ is equal to the multiplicative monotone convolution of the distributions of $U$ and $V$, i.e.\
\begin{equation}\label{eq-mon-mult-conv}
\mathcal{L}(UV,\Omega) = \mathcal{L}(VU,\Omega)  = \mathcal{L}(U,\Omega)\mmultg\mathcal{L}(V,\Omega).
\end{equation}
\end{theorem}
\begin{remark}
Note that the order of the convolution product on the right-hand-side of Equation \eqref{eq-mon-mult-conv} depends only on the order in which the operators $U-\mathbf{1}$ and $V-\mathbf{1}$ are monotonically independent, but not on the order in which $U$ and $V$ are multiplied.
\end{remark}

\section{Boolean Convolutions}

\begin{definition}
Let $\mathcal{A}_1,\mathcal{A}_2\subset\mathcal{B}(H)$ be two
$*$-algebras of bounded operators on a Hilbert space $H$, and let
$\Omega\in H$ be a unit vector. We say that $\mathcal{A}_1$ and
$\mathcal{A}_2$ are boolean independent w.r.t.\ $\Omega$, if we have
\[
\langle \Omega, X_1X_2\cdots X_k\Omega\rangle = \prod _{\kappa=1}^k
\langle \Omega, X_\kappa\Omega\rangle
\]
for all $k\in\mathbb{N}$, $\varepsilon\in\mathbb{A}_k$,
$X_1\in\mathcal{A}_{\varepsilon_1},\ldots,X_k\in\mathcal{A}_{\varepsilon_k}$.
\end{definition}

\begin{remark}
The algebras are not required to be unital. If one of them is
unital, say $\mathcal{A}_1$, then the restriction of
$\Phi(\cdot)=\langle\Omega,\,\cdot\,\Omega\rangle$ to the other
algebra, say $\mathcal{A}_2$, has to be a homomorphism, since the
boolean independence implies
\[
\langle \Omega, XY\Omega\rangle =  \langle \Omega,
X\mathbf{1}Y\Omega\rangle  = \langle \Omega, X\Omega\rangle\langle
\Omega,Y\Omega\rangle
\]
for $X,Y\in\mathcal{A}_2$.
\end{remark}

\begin{definition}\label{def-bool-indep}
Let $X$ and $Y$ be two normal operators on a Hilbert space $H$,
not necessarily bounded. We say that $X$ and $Y$ are boolean
independent, if the $*$-algebras ${\rm alg}_0(X)=\{h(X):h\in C_b(
\mathbb{C}),h(0)=0\}$ and ${\rm alg}_0(Y)=\{h(Y):h\in C_b(
\mathbb{C}),h(0)=0\}$ are boolean independent.
\end{definition}

We will start by characterizing up to unitary transformations the general form of two boolean independent normal operators. Given a measure space $(M,\mathcal{M},\mu)$, we shall denote by $L^2(M,\mu)_0$ the orthogonal complement of the constant function, i.e.\
\[
L^2(M,\mu)_0 = \left\{\psi\in L^2(M,\mu); \int_M \psi {\rm d}\mu=0\right\}.
\]

\begin{proposition}\label{prop-bool-model}
Let $\mu,\nu$ be two probability measures on $\mathbb{C}$ and define normal operators $N_x$ and $N_y$ on $\mathbb{C}\oplus L^2(\mathbb{C},\mu)_0\oplus L^2(\mathbb{C},\nu)_0$ by
\begin{eqnarray*}
{\rm Dom}\, N_x &=& \left\{\left(\begin{array}{c} \alpha \\ \psi_1 \\ \psi_2 \end{array}\right)\in \mathbb{C}\oplus L^2(\mathbb{C},\mu)_0\oplus L^2(\mathbb{C},\nu)_0; \int_{\mathbb{C}} \left| x \big(\psi_1(x)+\alpha\big)\right|^2{\rm d}\mu(x) < \infty\right\}, \\
{\rm Dom}\, N_y &=& \left\{\left(\begin{array}{c} \alpha \\ \psi_1 \\ \psi_2 \end{array}\right)\in \mathbb{C}\oplus L^2(\mathbb{C},\mu)_0\oplus L^2(\mathbb{C},\nu)_0; \int_{\mathbb{C}} \left| y \big(\psi_2(y)+\alpha\big)\right|^2{\rm d}\nu(y) < \infty\right\},
\end{eqnarray*}
\begin{eqnarray*}
N_x\left(\begin{array}{c} \alpha \\ \psi_1 \\ \psi_2 \end{array}\right) &=& \left(\begin{array}{c} \int_{\mathbb{C}}x\big(\psi_1(x)+\alpha\big){\rm d}\mu(x) \\ x(\psi_1+\alpha)- \int_{\mathbb{C}}x\big(\psi_1(x)+\alpha\big){\rm d}\mu(x) \\ 0 \end{array}\right), \\[1mm]
N_y\left(\begin{array}{c} \alpha \\ \psi_1 \\ \psi_2 \end{array}\right)&=& \left(\begin{array}{c} \int_{\mathbb{C}}y\big(\psi_2(y)+\alpha\big){\rm d}\nu(y) \\ 0 \\ x(\psi_2+\alpha) - \int_{\mathbb{C}}y\big(\psi_2(y)+\alpha\big){\rm d}\nu(y) \end{array}\right).
\end{eqnarray*}
Then $N_x$ and $N_y$ are boolean independent w.r.t.\ the vector $\omega=\left(\begin{array}{c} 1 \\ 0 \\ 0 \end{array}\right)$ and we have $\mathcal{L}(N_x,\omega) = \mu$, $\mathcal{L}(N_y,\omega) = \nu$.
\end{proposition}
\begin{proof}
Under the identification $\mathbb{C}\oplus L^2(\mathbb{C},\mu)_0\oplus L^2(\mathbb{C},\nu)_0\cong L^2(\mathbb{C},\mu)\oplus L^2(\mathbb{C},\nu)_0$, where
\[
\left(\begin{array}{c} \alpha \\ \psi_1 \\ \psi_2 \end{array}\right) \cong \left(\begin{array}{c} \psi_1+\alpha \\ \psi_2 \end{array}\right),
\]
the operator $N_x$ becomes multiplication by the variable $x$ on $L^2(\mathbb{C},\mu)$. It is clearly normal and we have
\[
h(N_x)\left(\begin{array}{c} \alpha \\ \psi_1 \\ \psi_2
\end{array}\right) =
\left( \begin{array}{c} \int_{\mathbb{C}} h(x)\big(\alpha + \psi_1(x)\big){\rm d}\mu(x) \\
h(\alpha+\psi_1) -\int_{\mathbb{C}} h(x)\big(\alpha + \psi_1(x)\big){\rm d}\mu(x) \\
h(0)\psi_2   \end{array}\right)
\]
and $\langle \omega, h(N_x)\omega\rangle=\int_{\mathbb{C}} h{\rm d}\mu$ for all $h\in C_b(\mathbb{C})$, i.e.\ $\mathcal{L}(N_x,\omega)=\mu$. Similarly
\[
h(N_y)\left(\begin{array}{c} \alpha \\ \psi_1 \\ \psi_2
\end{array}\right) =   \left( \begin{array}{c}\int_{\mathbb{C}} h(y)\big(\alpha + \psi_2(y)\big){\rm d}\nu(y) \\
h(0)\psi_1 \\
h(\alpha+\psi_2) -\int_{\mathbb{C}} h(y)\big(\alpha + \psi_2(y)\big){\rm d}\nu(y)   \end{array}\right)
\]
for all $h\in C_b(\mathbb{C})$, and $\mathcal{L}(N_y,\omega)=\nu$.

Let $f_1,\ldots,f_n,g_1,\ldots,g_n\in C_b(\mathbb{C})$, with $f_1(0)=\cdots=f_n(0)=g_1(0)=\cdots=g_n(0)=0$. Then
\[
f_n(N_x)g_{n-1}(N_y)\cdots g_1(N_y)f_1(N_x)\omega = \left(\begin{array}{c} \prod_{k=1}^n \int_{\mathbb{C}} f_k{\rm d}\mu \prod_{\ell=1}^{n-1} \int_{\mathbb{C}} g_\ell{\rm d}\nu 
\\
\prod_{k=1}^{n-1} \int_{\mathbb{C}} f_k{\rm d}\mu \prod_{\ell=1}^{n-1} \int_{\mathbb{C}} g_\ell{\rm d}\nu \left(f_n - \int_{\mathbb{C}}f_n{\rm d}\mu\right) \\
0\end{array}\right)
\]
and therefore
\begin{eqnarray*}
\langle \omega, f_n(N_x)g_{n-1}(N_y)\cdots g_1(N_y)f_1(N_x)\omega\rangle 
&=& \prod_{k=1}^n \int_{\mathbb{C}} f_k\,{\rm d}\mu \,\prod_{\ell=1}^{n-1} \int_{\mathbb{C}} g_\ell\,{\rm d}\nu  \\
&=& \prod_{k=1}^n \langle\omega,f_k(N_x)\omega\rangle \prod_{\ell=1}^{n-1} \langle\omega,g_\ell(N_y)\omega\rangle
\end{eqnarray*}
i.e.\ the condition for boolean independence is satisfied in this case.
Similarly one checks the expectation of $g_n(N_y)f_n(N_x)\cdots g_1(N_y)f_1(N_x)$, $f_n(N_x)g_n(N_y)\cdots f_1(N_x)g_1(N_y)$, and $g_n(N_y)f_{n-1}(N_x)\cdots f_1(N_x)g_1(N_y)$.
\end{proof}

We shall now show that any pair of boolean independent normal operators can be reduced to this model.

\begin{theorem}\label{thm-bool-model}
Let $X$ and $Y$ be two normal operators on a Hilbert space
$H$ that are boolean independent w.r.t.\ to $\Omega\in H$ and let
$\mu=\mathcal{L}(X,\Omega)$, $\nu=\mathcal{L}(Y,\Omega)$.

Then there exists an isometry $W:\mathbb{C}\oplus
L^2(\mathbb{C},\mu)_0\oplus L^2(\mathbb{C},\nu)_0 \to H$ such that
\begin{eqnarray}\label{eq-bool-model}
\\
W^*h(X)W\left(\begin{array}{c} \alpha \\ \psi_1 \\ \psi_2
\end{array}\right) &=&
\left( \begin{array}{c} \int_{\mathbb{C}} h(x)\big(\alpha + \psi_1(x)\big){\rm d}\mu(x) \\
h(\alpha+\psi_1) -\int_{\mathbb{C}} h(x)\big(\alpha + \psi_1(x)\big){\rm d}\mu(x) \\
h(0)\psi_2   \end{array}\right), \nonumber \\
W^*h(Y)W   \left(\begin{array}{c} \alpha \\ \psi_1 \\ \psi_2
\end{array}\right) &=&   \left( \begin{array}{c}\int_{\mathbb{C}} h(y)\big(\alpha + \psi_2(y)\big){\rm d}\nu(y) \\
h(0)\psi_1 \\
h(\alpha+\psi_2) -\int_{\mathbb{C}} h(y)\big(\alpha + \psi_2(y)\big){\rm d}\nu(y)   \end{array}\right) \nonumber
\end{eqnarray}
for all $h\in C_b(\mathbb{C})$, $\alpha\in\mathbb{C}$, $\psi_1\in
L^2(\mathbb{C},\mu)_0$, $\psi_2\in L^2(\mathbb{C},\nu)_0$.

We have  $W\big(\mathbb{C}\oplus L^2(\mathbb{C},\mu)_0\oplus L^2(\mathbb{C},\nu)_0\big)=\overline{{\rm alg}\{h(X),h(Y):h\in C_b( \mathbb{C})\} \Omega}$.

If the vector $\Omega\in H$ is cyclic for the algebra ${\rm alg}(X,Y)={\rm alg}\{h(X),h(Y):h\in C_b( \mathbb{C})\}$ generated by $X$ and $Y$, then $W$ is unitary.
\end{theorem}
\begin{proof}
For a probability measure $\mu$ on $\mathbb{C}$, let
\[
C_b(\mathbb{C})_{\mu,0} = \left\{ f\in C_b(\mathbb{C});
\int_\mathbb{C} f(z){\rm d}\mu(x) =0\right\},
\]
then $C_b(\mathbb{C})_{\mu,0}$ is dense in $L^2(\mathbb{C},\mu)_0$.

Define $W:\mathbb{C}\oplus C_b(\mathbb{C})_{\mu,0}\oplus
C_b(\mathbb{C})_{\nu,0}\to H$ by
\[
W \left(\begin{array}{c}\alpha \\ f \\ g\end{array}\right) =
\big(\alpha+f(X)+g(Y)\big)\Omega.
\]
This is an isometry, since
\begin{eqnarray*}
\left\langle W \left(\begin{array}{c}\alpha_1 \\ f_1 \\
g_1\end{array}\right),
W \left(\begin{array}{c}\alpha_2 \\ f_2 \\
g_2\end{array}\right)\right\rangle &=&
\left\langle\big(\alpha_1+f_1(X)+g_1(Y)\big)\Omega,
\big(\alpha_2+f_2(X)+g_2(Y)\big)\Omega\right\rangle \\
&=&\overline{\alpha_1}\alpha_2 + \int_{\mathbb{C}}
\overline{f_1(x)}f_2(x){\rm d}\mu(x) + \int_{\mathbb{C}}
\overline{g_1(y)}g_2(y){\rm d}\mu(y),
\end{eqnarray*}
where the mixed terms all vanish because
$\langle\Omega,f_i(X)\Omega\rangle=\langle\Omega,g_i(Y)\Omega\rangle
= 0$ for $i=1,2$. Therefore $W$ extends in a unique way to an
isometry on $\mathbb{C}\oplus L^2(\mathbb{C},\mu)_0\oplus
L^2(\mathbb{C},\nu)_0$

Let now $h\in C_b(\mathbb{C})$, then we get
\begin{multline*}
\left\langle W \left(\begin{array}{c}\alpha_1 \\ f_1 \\
g_1\end{array}\right),
h(X)W \left(\begin{array}{c}\alpha_2 \\ f_2 \\
g_2\end{array}\right)\right\rangle \\
=
\left\langle\big(\alpha_1+f_1(X)+g_1(Y)\big)\Omega,
(h(X)-h(0)\mathbf{1}\big) \big(\alpha_2+f_2(X)+g_2(Y)\big)\Omega\right\rangle \\
+ h(0)\left\langle W \left(\begin{array}{c}\alpha_1 \\ f_1 \\
g_1\end{array}\right),W \left(\begin{array}{c}\alpha_2 \\ f_2 \\
g_2\end{array}\right)\right\rangle\\
= \left\langle \big(\alpha_1+f_1(X)\big)\Omega, \big(h(X)-h(0)\mathbf{1}\big)\big(\alpha_2+f_2(X)\big)\Omega\right\rangle +h(0)\left\langle \left(\begin{array}{c}\alpha_1 \\ f_1 \\
g_1\end{array}\right),\left(\begin{array}{c}\alpha_2 \\ f_2 \\
g_2\end{array}\right)\right\rangle,
\end{multline*}
because the boolean independence and $\langle\Omega,g_i(Y)\Omega\rangle=0$ imply that all other terms vanish. But since $\langle\Omega,f_i(Y)\Omega\rangle=0$, this is equal to
\begin{multline*}
 \left\langle \big(\alpha_1+f_1(X)\big)\Omega, h(X)\big(\alpha_2+f_2(X)\big)\Omega\right\rangle \\
+h(0)\left(\left\langle \left(\begin{array}{c}\alpha_1 \\ f_1 \\
g_1\end{array}\right),\left(\begin{array}{c}\alpha_2 \\ f_2 \\
g_2\end{array}\right)\right\rangle - \overline{\alpha_1}\alpha_2-\langle f_1,f_2\rangle\right) \\
=\left\langle\left(\begin{array}{c} \alpha_1 \\ f_1 \\g_1 \end{array}\right) , \left(\begin{array}{c} \int h(x)\big(f_2(x)+\alpha_2\big){\rm d}\mu(x) \\ h(f_2+\alpha_2) - \int h(x)\big(f_2(x)+\alpha_2\big){\rm d}\mu(x) \\ h(0)g_2 \end{array}\right)\right\rangle.
\end{multline*}
This proves the first formula. The second formula follows by symmetry.

Let $f,g\in C_b(\mathbb{C})$, $f(0)=0$, and note that 
\begin{eqnarray*}
&& \left|\left|f(X)g(Y)\Omega - \int_{\mathbb{C}}g{\rm d}\nu\,f(X)\Omega\right|\right|^2 \\
&=& \langle \Omega, g(Y)^*|f(X)|^2g(Y)\Omega - \int_{\mathbb{C}}g{\rm d}\nu\,\langle\Omega, g(Y)^*|f(X)|^2 \Omega \\
&& - \int_{\mathbb{C}}\overline{g}{\rm d}\nu\,\langle\Omega, |f(X)|^2g(Y) \Omega\rangle +  \left(\int_{\mathbb{C}}g{\rm d}\nu\right)^2 \langle\Omega,|f(X)|^2\Omega\rangle \\
&=& 0,
\end{eqnarray*}
i.e.\ $f(X)g(Y)\Omega= \int_{\mathbb{C}}g{\rm d}\nu\,f(X)\Omega$. Similarly $f(Y)g(X)\Omega = \int_{\mathbb{C}}g{\rm d}\mu\, f(Y)\Omega$
and thus
\begin{eqnarray*}
\overline{{\rm alg}\{h(X),h(Y):h\in C_b( \mathbb{C})\} \Omega} &=& \overline{{\rm span}\,\{\Omega,f(X)\Omega,f(Y)\Omega; f\in C_b(\mathbb{C})\}} \\
&=& W\big(\mathbb{C}\oplus L^2(\mathbb{C},\mu)_0\oplus L^2(\mathbb{C},\nu)_0\big).
\end{eqnarray*}
If $\Omega$ is cyclic, then $W$ is surjective and therefore unitary.
\end{proof}
\begin{remark}\label{rem-bool-reduction}
As in the monotone case, cf.\ Remark \ref{rem-mon-reduction}, this theorem shows that joint law of bounded functions on $X$ and $Y$ is uniquely determined by $\mathcal{L}(X,\Omega)$ and $\mathcal{L}(Y,\Omega)$. Furthermore, the characterisation and computation of the law of unbounded functions of $X$ and $Y$ like, e.g., $X+Y$ or $\sqrt{X}Y\sqrt{Y}$, is also reduced to the model introduced in Proposition \ref{prop-bool-model}.
\end{remark}

\subsection{Additive boolean convolution on
$\mathcal{M}_1(\mathbb{R})$}

\begin{definition}\cite{speicher+woroudi93}
Let $\mu$ and $\nu$ be two probability measures on $\mathbb{R}$ with
reciprocal Cauchy transforms $F_\mu$ and $F_\nu$. Then we define the
additive monotone convolution $\lambda=\mu\uplus\nu$ of $\mu$ and
$\nu$ as the unique probability measure $\lambda$ on $\mathbb{R}$
with reciprocal Cauchy transform given by
\[
F_\lambda(z)= F_\mu(z)+F_\nu(z)-z
\]
for $z\in\mathbb{C}^+$.
\end{definition}
That the additive boolean convolution is well-defined follows from
Subsection \ref{nevanlinna}. It is commutative and associative, $*$-weakly continuous, but not affine, cf.\ \cite{speicher+woroudi93}.

\begin{proposition}\label{prop-bool-add-resolvent}
Let $\mu$ and $\nu$ be two probabilities on $\mathbb{R}$ and define operators $N_x$ and $N_y$ as in Proposition \ref{prop-bool-model}. Then $N_x$ and $N_y$ are self-adjoint and boolean independent w.r.t.\ $\omega=\left(\begin{array}{c} 1 \\ 0 \\ 0\end{array}\right)$. Furthermore, the operator $z-N_x-N_y$ has a bounded inverse for all $z\in\mathbb{C}\backslash\mathbb{R}$, given by
\begin{equation}\label{eq-bool-add-resolvent}
(z-N_x-N_y)^{-1} \left(\begin{array}{c} \alpha \\ \psi_1 \\ \psi_2\end{array}\right) =\left(\begin{array}{c} \beta \\ \frac{\psi_1 + \beta x - c_x}{z-x} \\\frac{\psi_2 + \beta y - c_y}{z-y} \end{array}\right),
\end{equation}
where
\begin{equation}\label{eq-beta-value}
\beta=\frac{\alpha G_\mu(z)G_\nu(z) + G_\nu(z)\int_\mathbb{R}\frac{\psi_1(x)}{z-x}{\rm d}\mu(x) + G_\mu(z)\int_\mathbb{R} \frac{\psi_2(y)}{z-y}{\rm d}\nu(y)}{G_\mu(z)+G_\nu(z)-zG_\mu(z)G_\nu(z)},
\end{equation}
and $c_x,c_y\in\mathbb{C}$ have to be chosen such that
\begin{equation}\label{eq-cx-cy-constants}
\int_\mathbb{R}\frac{\psi_1(x) + \beta x - c_x}{z-x}{\rm d}\mu(x)=0=\int_\mathbb{R}\frac{\psi_2(y) + \beta y - c_y}{z-y}{\rm d}\nu(y).
\end{equation}
\end{proposition}
Note that Equation \eqref{eq-cx-cy-constants} yields the following formulas for the constants $c_x,c_y$,
\begin{eqnarray*}
c_x &=& \frac{\int\frac{\psi_1(x)}{z-x}{\rm d}\mu(x)+\beta\big(zG_\mu(z)-1\big)}{G_\mu(z)}, \\
c_y &=& \frac{\int\frac{\psi_2(y)}{z-y}{\rm d}\nu(y)+\beta\big(zG_\nu(z)-1\big)}{G_\nu(z)}.
\end{eqnarray*}
\begin{proof}
$N_x$ and $N_y$ are boolean independent by Proposition \ref{prop-bool-model}.

For $z\in\mathbb{C}^+$, we have ${\rm Im}\,F_\mu(z)\ge {\rm Im}\,z>0$, ${\rm Im}\,F_\nu(z)\ge {\rm Im}\,z>0$, and therefore
\[
{\rm Im}\frac{G_\mu(z)+G_\nu(z)-zG_\mu(z)G_\nu(z)}{G_\mu(z)G_\nu(z)} = {\rm Im}\,\big(F_\mu(z)+F_\nu(z)-z\big)>0.
\]
This shows that the denominator of the right-hand-side of Equation \eqref{eq-beta-value} can not vanish for $z\in\mathbb{C}^+$. Since $G_\mu(\overline{z}) = \overline{G_\mu(z)}$, $G_\nu(\overline{z}) = \overline{G_\nu(z)}$, it can not vanish for $z$ with ${\rm Im}\,z<0$, either. The functions $\frac{1}{z-x}$ and $\frac{x}{z-x}$ are bounded on $\mathbb{R}$ for $z\in\mathbb{C}\backslash\mathbb{R}$, therefore Equation \eqref{eq-bool-add-resolvent} defines a bounded operator.

Let
\[
\varphi_1=\frac{\psi_1 + \beta x - c_x}{z-x} \qquad\mbox{ and }\qquad \varphi_2=\frac{\psi_2 + \beta y - c_y}{z-y},
\]
then
\[
(z-N_x-N_y)\left(\begin{array}{c} \beta \\ \varphi_1 \\ \varphi_2 \end{array}\right) = \left( \begin{array}{c} z\beta + d_x+d_y \\ (z-x)\varphi_1 - \beta x - d_x \\ (z-y)\varphi_2 - \beta y - d_y\end{array}\right) =\left( \begin{array}{c} z\beta + d_x+d_y \\ \psi_1 - c_x - d_x \\ \psi_2 - c_y - d_y\end{array}\right) \\
\]
where
\[
d_x= \int x\big(\varphi_1(x)+\beta\big){\rm d}\mu(x), \qquad d_y= \int y \big(\varphi_2(y)+\beta\big){\rm d}\nu(y).
\]
Since $\psi_1\in L^2(\mathbb{R},\mu)_0$, $\psi_2\in L^2(\mathbb{R},\nu)_0$, integrating over the second and third component gives $c_x=-d_x$ and $c_y=-d_y$. Therefore
\[
(z-N_x-N_y)\left(\begin{array}{c} \beta \\ \varphi_1 \\ \varphi_2 \end{array}\right) = \left( \begin{array}{c} z\beta - c_x-c_y \\ \psi_1 \\ \psi_2 \end{array}\right) \\
\]
We have to show that the first component is equal to $\alpha$. We get
\begin{multline*}
z\beta - c_x-c_y = z\beta -  \frac{\int\frac{\psi_1(x)}{z-x}{\rm d}\mu(x)+\beta\big(zG_\mu(z)-1\big)}{G_\mu(z)} -  \frac{\int\frac{\psi_1(x)}{z-x}{\rm d}\mu(x)+\beta\big(zG_\nu(z)-1\big)}{G_\nu(z)} \\
= \beta\frac{G_\mu(z)+G_\nu(z)-zG_\mu(z)G_\nu(z)}{G_\mu(z)G_\nu(z)} - \frac{1}{G_\mu(z)}\int\frac{\psi_1(x)}{z-x}{\rm d}\mu(x) - \frac{1}{G_\nu(z)}\int\frac{\psi_2(y)}{z-y}{\rm d}\nu(y)
\end{multline*}
Substituting Equation \eqref{eq-beta-value} into this expression, we get the desired result $z\beta-c_x-c_y=\alpha$.
\end{proof}
\begin{theorem}\label{thm-bool-add-conv}
Let $X$ and $Y$ be two self-adjoint operators on a Hilbert space $H$ that are
boolean independent w.r.t.\ a unit vector $\Omega\in H$ and assume that $\Omega$
is cyclic, i.e.\ that
\[
\overline{{\rm alg}\{h(X),h(Y);h\in C_b(\mathbb{R})\}\Omega}=H.
\]
Then $X+Y$ is essentially self-adjoint and the distribution w.r.t.\ $\Omega$ of the closure of $X+Y$ is equal to the boolean convolution of the distributions of $X$ and $Y$ w.r.t.\ $\Omega$, i.e.\
\[
\mathcal{L}(X+Y,\Omega) = \mathcal{L}(X,\Omega)\uplus\mathcal{L}(Y,\Omega).
\]
\end{theorem}
\begin{proof}
Let $\mu=\mathcal{L}(X,\Omega)$, $\nu=\mathcal{L}(Y,\Omega)$.

By Theorem \ref{thm-bool-model} and Lemma \ref{lem-unitary-eq} it is sufficient to consider the case where $X$ and $Y$ are defined as in Proposition \ref{prop-bool-model}. Then Proposition \ref{prop-bool-add-resolvent} shows that $z-X-Y$ admits a bounded inverse for all $z\in\mathbb{C}\backslash\mathbb{R}$ and therefore that ${\rm Ran}\,(z-X-Y)$ is dense. By \cite[Theorem VIII.3]{reed+simon80} this is equivalent to $X+Y$ being essentially self-adjoint.

Using Equation \eqref{eq-bool-add-resolvent}, we can compute the Cauchy transform of the distribution of the closure of $X+Y$. Let $z\in\mathbb{C}^+$, then
\begin{eqnarray*}
G_{X+Y}(z) &=& \langle\Omega, (z-X-Y)^{-1}\Omega\rangle = \left\langle \omega, (z-N_x-N_y)^{-1}\omega\right\rangle \\
&=&  \left\langle\left(\begin{array}{c}1 \\ 0 \\ 0  \end{array}\right), \frac{G_\mu(z)G_\nu(z)}{G_\mu(z)+G_\nu(z)-zG_\mu(z)G_\nu(z)}\left(\begin{array}{c} 1 \\ \frac{x-\frac{zG_\mu(z)-1}{G_\mu(z)}}{z-x} \\ \frac{y-\frac{zG_\nu(z)-1}{G_\nu(z)}}{z-y} \end{array}\right)\right\rangle \\
&=& \frac{G_\mu(z)G_\nu(z)}{G_\mu(z)+G_\nu(z)-zG_\mu(z)G_\nu(z)}.
\end{eqnarray*}
Replacing all Cauchy transforms by their reciprocals, this becomes
\[
F_{X+Y}(z) = F_\mu(z)+F_\nu(z)-z=F_{\mu\uplus\nu}(z).
\]
\end{proof}

\subsection{Multiplicative boolean convolution on
$\mathcal{M}_1(\mathbb{R}_+)$}

Let us first recall Bercovici's definition of the boolean convolution for probability measures in the positive half-line, cf.\ \cite{bercovici04b}.

\begin{definition}\cite{bercovici04b}
Let $\mu$ and $\nu$ be two probability measures on $\mathbb{R}_+$ with transforms $K_\mu$ and $K_\nu$. If the holomorphic function defined by
\begin{equation}\label{eq-def-mult1-bool}
K(z) = \frac{K_\mu(z)K_\nu(z)}{z}
\end{equation}
for $z\in\mathbb{C}\backslash\mathbb{R}_+$ belongs to the class $\mathcal{P}$ introduced in Subsection \ref{nevanlinna}, then the boolean convolution $\lambda=\mu\bmultg\nu$ is defined as the unique probability measure $\lambda$ on $\mathbb{R}_+$ with transform $K_\lambda=K$.
\end{definition}
But in general the function $K$ defined in Equation \eqref{eq-def-mult1-bool} does not belong to $\mathcal{P}$ and in that case the convolution of $\mu$ and $\nu$ is not defined. Bercovici has shown that for any probability measure $\mu$ on $\mathbb{R}_+$ not concentrated in one point there exists an $n\in\mathbb{N}$ such that the $n$-fold convolution product $\mu^{\bmultk n}$ of $\mu$ with itself is not defined, cf.\ \cite[Proposition 3.1]{bercovici04b}.

This is of course related to the problem that in general the product of two positive operators is not positive. One might hope that taking e.g.\ $\sqrt{X}Y\sqrt{X}$ could lead to a better definition of the multiplicative boolean convolution, since this operator will automatically be positive. We will see below that this approach leads to a convolution that is always defined, but that is not associative.

Let us first introduce the model which we will use for our calculations.

\begin{proposition}
Let $\mu$ and $\nu$ be probability measures on $\mathbb{R}_+$. Define operators $Q_x$ and $Q_y$ on $H=\mathbb{C}\oplus L^2(\mathbb{R}_+,\mu)_0\oplus L^2(\mathbb{R},\nu)_0$ by
\begin{eqnarray*}
{\rm Dom}\,Q_x &=&\left\{\left(\begin{array}{c} \alpha \\ \psi_1 \\ \psi_2 \end{array}\right)\in H ; \int_{\mathbb{R}_+} x^2\big|\psi_1(x)+\alpha\big|^2{\rm d}\mu(x) <\infty\right\}, \\
{\rm Dom}\,Q_y &=&\left\{\left(\begin{array}{c} \alpha \\ \psi_1 \\ \psi_2 \end{array}\right)\in H ; \int_{\mathbb{R}_+} y^2\big|\psi_2(y)+\alpha\big|^2{\rm d}\mu(y) <\infty\right\}, \\
Q_x\left(\begin{array}{c} \alpha \\ \psi_1 \\ \psi_2 \end{array}\right) &=& \left(\begin{array}{c} \int_{\mathbb{R}_+}x\big(\psi_1(x)+\alpha\big){\rm d}\mu(x) \\ x(\psi_1+\alpha) - \int_{\mathbb{R}_+}x\big(\psi_1(x)+\alpha\big){\rm d}\mu(x) \\ \psi_2 \end{array}\right), \\
Q_x\left(\begin{array}{c} \alpha \\ \psi_1 \\ \psi_2 \end{array}\right) &=& \left(\begin{array}{c} \int_{\mathbb{R}_+}y\big(\psi_2(y)+\alpha\big){\rm d}\nu(y) \\ \psi_1 \\ y(\psi_2+\alpha) - \int_{\mathbb{R}_+}y\big(\psi_2(y)+\alpha\big){\rm d}\nu(y) \end{array}\right).
\end{eqnarray*}
Then $Q_x-\mathbf{1}$ and $Q_y-\mathbf{1}$ are boolean independent w.r.t.\ $\omega=\left(\begin{array}{c} 1 \\ 0 \\ 0 \end{array}\right)$, and $\mathcal{L}(Q_x,\omega)=\mu$, $\mathcal{L}(Q_y,\omega)=\nu$.
\end{proposition}
\begin{proof}
It follows from Propisition \ref{prop-bool-model}, that $Q_x-\mathbf{1}$ and $Q_y-\mathbf{1}$ are boolean independent.

Note that the functional calculus for $Q_x$ is given by $C_b(\mathbb{R}_+)\ni f\mapsto f(Q_x)\in\mathcal{B}(H)$,
\[
f(Q_x)\left(\begin{array}{c} \alpha \\ \psi_1 \\ \psi_2 \end{array}\right) = \left(\begin{array}{c} \int_{\mathbb{R}_+} f(x)\big(\psi_1(x)+\alpha\big){\rm d}\mu(x) \\ f(\psi_1+\alpha) -\int_{\mathbb{R}_+} f(x)\big(\psi_1(x)+\alpha\big){\rm d}\mu(x) \\ f(1)\psi_2 \end{array}\right).
\]
Therefore
\[
\langle \omega , f(Q_x)\omega\rangle = \left\langle\left(\begin{array}{c} 1 \\ 0 \\ 0 \end{array}\right), \left(\begin{array}{c} \int_{\mathbb{R}_+}f(x){\rm d}\mu(x) \\ f-\int_{\mathbb{R}_+}f(x){\rm d}\mu(x) \\ 0 \end{array}\right)\right\rangle = \int_{\mathbb{R}_+}f(x){\rm d}\mu(x),
\]
i.e.\ $\mathcal{L}(X,\omega)=\mu$. Similarly $\mathcal{L}(Y,\omega)=\nu$.
\end{proof}

\begin{proposition}
Let $z\in\mathbb{C}\backslash\mathbb{R}$, then $z-\sqrt{X}Y\sqrt{X}$ has a bounded inverse, given by
\begin{equation}\label{eq-bool-mult-xyx}
\left(z-\sqrt{X}Y\sqrt{X}\right)^{-1}\left(\begin{array}{c}\alpha \\ \psi_1 \\ \psi_2 \end{array}\right) = \left(\begin{array}{c} \beta \\ \varphi_1 \\ \varphi_2 \end{array}\right),
\end{equation}
where
\begin{eqnarray*}
\varphi_1(x) &=& \frac{\psi_1(x) + \beta x + (c_2-c_1)\sqrt{x}-c_3}{z-x}, \\
\varphi_2(y) &=& \frac{\psi_2(y) +c_1y-c_2}{z-y}, \\
\beta &=& \alpha G_\mu(z) + \int_{\mathbb{R}_+}\frac{\psi_1(x)}{z-x}{\rm d}\mu(x) + W_\mu(z)(c_2-c_1), \\
c_1 &=& \frac{\big(zG_\mu(z)-1\big)\int_{\mathbb{R}_+}\frac{\psi_2(y)}{z-y}{\rm d}\nu(y) + G_\nu(z)\left(\alpha W_\mu(z) + \int_{\mathbb{R}_+}\frac{\sqrt{x}\psi_1(x)}{z-x}{\rm d}\mu(x)\right)}{zG_\mu(z)G_\nu(z)-\big(zG_\mu(z)-1\big)\big(zG_\nu(z)-1\big)}, \\
c_2 &=&  \frac{zG_\mu(z)\int_{\mathbb{R}_+}\frac{\psi_2(y)}{z-y}{\rm d}\nu(y) + \big(zG_\nu(z)-1\big)\left(\alpha W_\mu(z) + \int_{\mathbb{R}_+}\frac{\sqrt{x}\psi_1(x)}{z-x}{\rm d}\mu(x)\right)}{zG_\mu(z)G_\nu(z)-\big(zG_\mu(z)-1\big)\big(zG_\nu(z)-1\big)}, \\
c_3 &=& z\beta-\alpha,
\end{eqnarray*}
and $W_\mu$ denotes the Cauchy transform of $\sqrt{x}\mu$, i.e.\
\[
W_\mu(z) = G_{\sqrt{x}\mu}(z) = \int_{\mathbb{R}_+}\frac{\sqrt{x}}{z-x}{\rm d}\mu(x).
\]
\end{proposition}
\begin{proof}
This can be checked by applying $z-\sqrt{X}Y\sqrt{X}$ to the right-hand-side of Equation \eqref{eq-bool-mult-xyx}. The computations are straight-forward, but rather tedious.
\end{proof}

\begin{theorem}
Let $X$ and $Y$ be two positive operators on a Hilbert space $H$ such that $X-\mathbf{1}$ and $Y-\mathbf{1}$ are boolean independent w.r.t.\ a unit vector $\Omega\in H$. Suppose furthermore that $\Omega$ is cyclic, i.e.\ that
\[
\overline{{\rm alg}\{h(X),h(Y);h\in C_b(\mathbb{R}_+)\}\Omega}=H.
\]
Then $\sqrt{X}Y\sqrt{X}$ is essentially self-adjoint, its closure is positive, and the distribution $\lambda=\mathcal{L}(\sqrt{X}Y\sqrt{X})$ w.r.t.\ $\Omega$ of its closure is given by
\begin{equation}\label{eq-def-bool-mult-conv1}
G_\lambda(z) = G_X(z) + \frac{\big(W_X(z)\big)^2\big((z-1)G_Y(z)-1\big)}{zG_X(z)G_Y(z)-\big(zG_X(z)-1\big)\big(zG_Y(z)-1\big)}
\end{equation}
for $z\in\mathbb{C}\backslash\mathbb{R}$, where
\[
W_X(z) = \left\langle \Omega, \frac{\sqrt{X}}{z-X}\Omega\right\rangle.
\]
\end{theorem}
\begin{proof}
As in the previous cases, the existence of a bounded inverse of $z-\sqrt{X}Y\sqrt{X}$ for $z\in\mathbb{C}\backslash\mathbb{R}$ implies that $\sqrt{X}Y\sqrt{X}$ is essentially self-adjoint. Furthermore it is clearly positive.

Using Equation \eqref{eq-bool-mult-xyx}, one can calculate the Cauchy transform of the distribution of $\sqrt{X}Y\sqrt{X}$.
\end{proof}
\begin{remark}\label{rem-new-bool-mult-conv}
One can now use Equation \eqref{eq-def-bool-mult-conv1} to define a ``quantum probabilistically motivated'' boolean convolution $\bmultsqrt$ for probability measures on $\mathbb{R}_+$. Let $\mu,\nu\in\mathcal{M}_1(\mathbb{R}_+)$, then $\lambda=\mu\bmultsqrt\nu$ is defined as the unique probability measure $\lambda$ on $\mathbb{R}_+$ such that
\[
G_\lambda(z) = G_\mu(z) + \frac{\big(W_\mu(z)\big)^2\big((z-1)G_\nu(z)-1\big)}{zG_\mu(z)G_\nu(z)-\big(zG_\mu(z)-1\big)\big(zG_\nu(z)-1\big)}
\]
for $z\in\mathbb{C}^+$.

This new convolution is defined for arbitrary probability measures $\mu$ and $\nu$ on $\mathbb{R}_+$, but it is neither associative nor commutative. For explicit calculations we use again the matrices $X$ and $Y$ introduced in Remark \ref{rem-mon-mult-alt}. The matrices $X-\mathbf{1}$ and $Y-\mathbf{1}$ are boolean independent w.r.t.\ to $\Omega$, see also Subsection \ref{subsec-rel-dirac}. Therefore
\begin{eqnarray*}
\big(p \delta_0 + (1-p)\delta_y\big)\bmultsqrt\delta_x &=& \mathcal{L}\big(\sqrt{Y}X\sqrt{Y},\Omega\big) \\
&=& \delta_x\mmultalt\big(p\delta_0+(1-p)\delta_y\big) \\
&=& p\delta_0 + (1-p)\delta_{y(xp+1-p)}, \\
\delta_x\bmultsqrt\big(p \delta_0 + (1-p)\delta_y\big) &=& \mathcal{L}\big(\sqrt{X}Y\sqrt{X},\Omega\big) \\
&=& \delta_x\mmultg\big(p\delta_0+(1-p)\delta_y\big) \\
&=& \frac{1-p}{xp+1-p}\delta_0 + \frac{xp}{xp+1-p}\delta_{y(xp+1-p)}
\end{eqnarray*}
for $x,y>0$, $0<p<1$. It is now easy to find explicit examples such that $\mu\bmultsqrt\nu\not=\nu\bmultsqrt\mu$ and $\lambda\bmultsqrt(\mu\bmultsqrt\nu)\not=(\lambda\bmultsqrt\mu)\bmultsqrt\nu$.
\end{remark}

\subsection{Multiplicative boolean convolution on $\mathcal{M}_1(\mathbb{T})$}\label{subsec-mult-bool-T}

For completeness we recall the results of \cite{franz04} for the multiplicative boolean convolution on $\mathcal{M}_1(\mathbb{T})$.

\begin{definition}\cite{franz04}
Let $\mu$ and $\nu$ be two probability measures on the unit circle $\mathbb{T}$ with transforms $K_\mu$ and $K_\nu$. Then the multiplicative monotone convolution $\lambda=\mu\bmultg\nu$ is defined as the unique probability on $\mathbb{T}$ with transform $K_\lambda$ given by
\[
K_\lambda(z) = \frac{K_\mu(z)K_\nu(z)}{z}
\]
for $z\in\mathbb{D}$.
\end{definition}

It is easy to deduce from Subsection \ref{nevanlinna} that the multiplicative boolean convolution on $\mathcal{M}_1(\mathbb{T})$ is well-defined. It is associative, commutative, $*$-weakly continuous in both arguments, but not affine.

\begin{theorem}\cite[Theorem 2.2]{franz04}
Let $U$ and $V$ be two unitary operators on a Hilbert space $H$, $\Omega\in H$ a unit vector and assume furthermore that $U-\mathbf{1}$ and $V-\mathbf{1}$ are boolean independent w.r.t.\ $\Omega$. Then the products $UV$ and $VU$ are also unitary and their distribution w.r.t. $\Omega$ is equal to the multiplicative boolean convolution of the distributions of $U$ and $V$, i.e.\
\[
\mathcal{L}(UV,\Omega) = \mathcal{L}(VU,\Omega)  = \mathcal{L}(U,\Omega)\bmultg\mathcal{L}(V,\Omega).
\]
\end{theorem}

\section{Some Relations between Free, Monotone, and Boolean Convolutions}\label{sec-relations}

\subsection{Decomposing free convolution products into monotone or boolean convolution products}

The theorems by Maassen, Chistyakov and G\"otze that we cited in Subsection \ref{free-conv} have an interesting formulation purely in terms of measures. We will consider only the additive case here, but similar results exist also for the two multiplicative free convolutions.

\begin{theorem}\label{thm-free-mon-bool-conv}
Let $\mu$ and $\nu$ be two probability measures on the real line. Then there exist two unique probability measures $\zeta_1$ and $\zeta_2$ on the real line such that
\[
\mu\boxplus\nu=\mu\triangleright\zeta_1 = \nu\triangleright\zeta_2 = \zeta_1\uplus\zeta_2.
\]
\end{theorem}
\begin{proof}
Apply Theorem \ref{thm-parallelogramm} to $F_1=F_\mu$ and $F_2=F_\nu$, and take for $\zeta_1$ and $\zeta_2$ the probability measures on $\mathbb{R}$ with reciprocal Cauchy transforms $Z_1$ and $Z_2$, respectively.
\end{proof}
\begin{remark}
The existence of unique probability measures $\zeta_1,\zeta_2$ such that $\mu\boxplus\nu=\mu\triangleright\zeta_1$ and $\mu\boxplus\nu=\nu\triangleright\zeta_2$ follows also from analytic subordination. E.g., $\zeta_1$ is obtained from the Markov kernel in \cite[Theorem 3.1]{biane98} by setting $x=0$.
\end{remark}
Recently, Accardi, Lenczewski, and Sa{\l}apata have given a similar result for products of graphs. Given two graphs $\mathcal{G}_1$ and $\mathcal{G}_2$, they gave an explicit construction of two more graphs $\mathcal{B}_1$ and $\mathcal{B}_2$ such that the free product of $\mathcal{G}_1$ and $\mathcal{G}_2$, the star product of $\mathcal{B}_1$ and $\mathcal{B}_2$, and the comb products of $\mathcal{G}_1$ and $\mathcal{B}_2$, or $\mathcal{G}_2$ and $\mathcal{B}_1$ are all isomorphic. For details, see \cite{accardi+lenczewski+salapata06}.

Lenczweski \cite{lenczewski06} has also given a more explicit version of Theorem \ref{thm-free-mon-bool-conv}. Given two bounded free operators $X$ and $Y$, he decomposes their sum as $X+Y=X_0+Z$ such that $X_0$ has the same distribution as $X$ and $X_0$ and $Z$ are monotonically independent. Let $H_X$ denote the subspace generated by ${\rm alg}(X)$ from the vacuum vector, $H_X=\overline{{\rm alg}(X)\Omega}$ and $P_X$ the orthogonal projection onto $H_X$. Then $X_0$ is given by $X_0=XP_X$ and $Z$ by $Z=X(\mathbf{1}-P_X)+Y$.

Theorem \ref{thm-free-mon-bool-conv} has also an interesting consequence for independent increment processes.
\begin{corollary}
Let $T>0$ and $(\mu_{st})_{0\le s\le t \le T}$ be a free convolution hemi-group in $\mathcal{M}_1(\mathbb{R})$, i.e.\ a two-parameter family of probability measures on the real line such that
\[
\mu_{st}\boxplus\mu_{tu}=\mu_{su}
\]
for all $0\le s\le t\le u\le T$. Then there exists a unique monotone convolution hemi-group $(\zeta_{st})_{0\le s\le t\le T}$ such that $\mu_{0t}=\zeta_{0t}$ for all $0\le t\le T$.
\end{corollary}
\begin{proof}
Let $0\le s\le t\le T$ and define $\zeta_{st}$ as the unique probability measure such that
\[
\mu_{0t}=\mu_{0s}\boxplus\mu_{st} = \mu_{0s}\triangleright\zeta_{st}.
\]
Clearly, we have $\zeta_{0t}=\mu_{0t}$ for all $0\le t\le T$. To check that the $\zeta_{st}$ form a monotone convolution hemi-group, rewrite $\mu_{0u}$ in two ways,
\begin{eqnarray*}
\mu_{0u} &=& \mu_{0s}\boxplus \mu_{su} = \mu_{0s}\triangleright\zeta_{su} \\
&=&
\mu_{0t}\boxplus\mu_{tu} = \mu_{0t}\triangleright\zeta_{tu} = (\mu_{0s}\boxplus\mu_{st})\triangleright\zeta_{tu} \\ &=& (\mu_{0s}\triangleright\zeta_{st})\triangleright\zeta_{tu} = \mu_{0s} \triangleright(\zeta_{st}\triangleright\zeta_{tu})
\end{eqnarray*}
and therefore by uniqueness $\zeta_{st}\triangleright\zeta_{tu}=\zeta_{su}$.
\end{proof}
Since independent increment processes are uniquely determined by the hemi-group of their marginal distributions, this induces a map from free independent increment processes to monotone independent increment processes. Under this map the free additive L\'evy processes of the second kind introduced in \cite{biane98} correspond exactly to stationary monotone increment processes.

\subsection{Monotone and boolean convolutions involving Dirac measures}\label{subsec-rel-dirac}

Let $X$ be a normal operator on some Hilbert space. The distribution of $X$ w.r.t.\ to some unit vector $\Omega\in H$ is concentrated in one point if and only if the vector state acts as a homomorphism on the algebra generated by $X$, or equivalently, if $\Omega$ is an eigenvector of $X$. Let $X$ and $Y$ be two normal operators and assume $\mathcal{L}(X,\Omega)=\delta_x$ for some $x\in\mathbb{C}$. Then $X$ and $Y$ are monotonically independent if and only if they are boolean independent. Therefore we get the following relations for monotone and boolean convolutions,
\begin{eqnarray*}
\delta_x\triangleright\mu &=& \delta_x\uplus\mu \qquad\mbox{ for } x\in\mathbb{R}, \quad \mu\in\mathcal{M}_1(\mathbb{R}), \\
\delta_x\mmultg\mu &=& \delta_x\bmultg\mu \qquad\mbox{ for } x\in\mathbb{T}, \quad \mu\in\mathcal{M}_1(\mathbb{T}), \\
\delta_x\mmultg\mu &=& \delta_x\bmultsqrt\mu \qquad\mbox{ for } x\in\mathbb{R}_+, \quad \mu\in\mathcal{M}_1(\mathbb{R}_+), \\
\delta_x\mmultalt\mu &=& \mu\bmultsqrt\delta_x \qquad\mbox{ for } x\in\mathbb{R}_+, \quad \mu\in\mathcal{M}_1(\mathbb{R}_+).
\end{eqnarray*}
{}From Equation \eqref{eq-int-formula} we now get, e.g.,
\[
\mu\triangleright\nu = \int_\mathbb{R} \delta_x\uplus\nu{\rm d}\mu(x)
\]
for $\mu,\nu\in\mathcal{M}_1(\mathbb{R})$, i.e.\ the monotone convolution can be considered as a linearization of the boolean convolution w.r.t.\ to the first argument.

\section*{Acknowledgements}

This work was completed while I was visiting the Graduate School of Information Sciences of Tohoku University as Marie-Curie fellow. I would like to thank Professors Nobuaki Obata, Fumio Hiai, and the other members of the GSIS for their hospitality. I am also indebted to an anonymous referee for suggesting important corrections and improvements.


\end{document}